\documentclass[12pt]{article}
\usepackage{amssymb,amsmath,amsthm}
\usepackage{amscd}
\usepackage{mathrsfs}
\newcommand{\N}{\mathbb{N}}
\newcommand{\R}{\mathbb{R}}

\newcommand{\dd}{\,{\rm d}}

\renewcommand{\div}{\mathop{\mathrm{div}}}

\numberwithin{equation}{section}
\newtheorem{thm}{Theorem}[section]

\newtheorem{prop}[thm]{Proposition}
\newtheorem{lem}[thm]{Lemma}
\newtheorem{rem}[thm]{Remark}

\begin{document}

\title{A Liouville theorem for the planer Navier-Stokes equations with the no-slip boundary condition and its application to a geometric regularity criterion}

\author{
\null\\
Yoshikazu Giga\\
Graduate School of Mathematical Sciences, University of Tokyo\\
3-8-1 Komaba, Meguro-ku, Tokyo 153-8914, Japan\\
{\tt labgiga@ms.u-tokyo.ac.jp}
\and
\\
Pen-Yuan Hsu\\
Graduate School of Mathematical Sciences, University of Tokyo\\
3-8-1 Komaba, Meguro-ku, Tokyo 153-8914, Japan\\
{\tt pyhsu@ms.u-tokyo.ac.jp}
\and
\\
Yasunori Maekawa\\
Mathematical Institute, Tohoku University\\          
6-3 Aoba, Aramaki, Aoba, Sendai 980-8578, Japan\\
{\tt maekawa@m.tohoku.ac.jp }
}

\date{}

\maketitle

\begin{abstract}
We establish a Liouville type result for a backward global solution to the Navier-Stokes equations in the half plane with the no-slip boundary condition.
No assumptions on spatial decay for the vorticity nor the velocity field are imposed. We study the vorticity equations instead of the original Navier-Stokes 
equations.  As an application, we extend the geometric regularity criterion for the Navier-Stokes equations in the three-dimensional half space under the no-slip boundary condition.
\end{abstract}

\section{Introduction}\label{sec.intro}
In this paper we study a backward solution to the Navier-Stokes equations in the half plane
\begin{align}
\partial_t u + {\div}\,  ( u \otimes u ) - \Delta u + \nabla p =0,~~~~~ {\rm div} \, u =0~~~~~~~~~{\rm in}~~ (-\infty, 0) \times \R^2_+\label{eq.NS.intro}
\end{align}
subject to the no-slip boundary condition
\begin{align}
u=0~~~~~~~~~~{\rm on} ~~ (-\infty, 0) \times \partial \R^2_+.\label{eq.noslip.intro}
\end{align}
Here $\R^2_+=\{(x_1,x_2)\in \R^2~|~ x_2>0\}$, and $u=u(t,x)=(u_1(t,x),u_2(t,x))$, $p=p(t,x)$ denote the velocity field, the pressure field, respectively.
We use the standard notation for derivatives;
$\partial_t = \partial/\partial t, ~\partial_j = \partial/\partial x_j,~\Delta = \sum_{j=1}^2 \partial_j^2,~{\div}\, u=\sum_{j=1}^2\partial_j u_j$,
and $(u \otimes u )_{1\leq i,j\leq 2} = (u_i u_j)_{1\leq i,j\leq 2}$.

We are interested in the Liouville problem for \eqref{eq.NS.intro} - \eqref{eq.noslip.intro}, that is, the nonexistence of nontrivial bounded global solutions to \eqref{eq.NS.intro} - \eqref{eq.noslip.intro}. 
As is well known, in the study of evolution equations the Liouville problem for bounded {\it backward} solutions plays an important role in obtaining an a priori bound of {\it forward} solutions
through a suitable scaling argument called a blow-up argument. For example, the reader is referred to \cite{Gi1}  for semilinear parabolic equations, to \cite{KoNaSeSv, SeSv1} 
for the axisymmetric Navier-Stokes equations  (see also \cite{ChStYaTs1, ChStYaTs2} for a different approach), 
to \cite{GiMi,Gi2} for a geometric regularity criterion to the three-dimensional Navier-Stokes equations, and to a recent result \cite{AbGi} for the Stokes semigroup in $L^\infty$ spaces.

This paper is particularly motivated by \cite{GiMi,Gi2}, where \eqref{eq.NS.intro} - \eqref{eq.noslip.intro} is naturally derived from a blow-up argument 
for the three-dimensional Navier-Stokes equations in the half space. 
Indeed,  if one imposes a uniform continuity on the alignment of the vorticity direction, the blow-up limit of the three-dimensional (Navier-Stokes) flow must be
a nontrivial bounded two-dimensional flow, and the problem is essentially reduced to the analysis of  \eqref{eq.NS.intro} - \eqref{eq.noslip.intro}.
If, in addition, one assumes that the possible blow-up is type I, then the limit flow is not allowed to be a constant in time. 
Thus the resolution of the Liouville problem is a crucial step to reach a contradiction.  From this systematic argument we can exclude the possibility of type I blow-up for the original three-dimensional flows 
under a regularity condition on the vorticity direction. 

Recently the paper \cite{GiMi} successfully completes the above argument when the velocity field satisfies the {\it perfect slip} boundary condition, but
the problem was remained open for the case of the {\it no-slip} boundary condition, which is physically more relevant. In this paper we prove a Liouville type theorem
for \eqref{eq.NS.intro} - \eqref{eq.noslip.intro} under some conditions on the velocity field $u$, the pressure field  $p$, and the vorticity field  $\omega=\partial_1 u_2-\partial_2u_1$. 
Our result is useful enough to settle the problem left open in \cite{GiMi}; see Theorem \ref{thm.intro.2} below. The details on this geometric regularity criterion will be discussed in Section \ref{sec.criterion}.

When one discusses the Liouville problem the choice of function spaces is of course a crucial issue. 
Indeed, if $u$ solves \eqref{eq.NS.intro} - \eqref{eq.noslip.intro} and decays fast enough in time and space then it is easy to conclude that $u$ is identically zero by a standard energy inequality. 
However, in view of application to the geometric regularity criterion, it is important to establish a Liouville type result within the framework of spatially nondecaying solutions. 
We should recall here that there are nontrivial shear flows whose velocity fields are bounded and decaying in time as $t\rightarrow -\infty$,
while the pressure fields grow linearly at spatial infinity; see \cite{SeSv2, Gi2}, and see also \eqref{eq.Poiseuille type flow.Sec.6} below.
The appearance of the  time-decaying shear flows is due to both the presence of the nontrivial boundary and the no-slip boundary condition in \eqref{eq.NS.intro} - \eqref{eq.noslip.intro}.
Indeed, if we consider the whole space case or if we replace \eqref{eq.noslip.intro} by the perfect slip boundary condition, $\partial_2u_1 = u_2=0$ on $\partial\R^2$, then such kind of flows does not exist. 
We note that these shear flows also solve the Stokes equations (i.e. nonlinear term is absent). 
Thus, even for the  linearized problem, we need to impose some assumptions on the spatial growth of the pressure field to obtain a Liouville theorem.
In fact, for the Stokes equations it is recently shown in \cite{JiSeSv} that any nontrivial bounded backward solution has to be a shear flow. 
Especially, the result of \cite{JiSeSv} gives a complete characterization of bounded backward solutions for the linear problem.

On the other hand, for the full Navier-Stokes equations there seems to be still few results on the Liouville type problem even in the case of the half plane.
The crucial difficulty is that, though the vorticity field satisfies the heat-transport equations,  
maximum principle  is no longer a useful tool to obtain an a priori bound of the vorticity field.
Indeed, the no-slip boundary condition on the velocity field is in general a source of vorticity on the boundary,
and maximum principle does not provide useful information about this vorticity production on the boundary.
This is contrasting with the case of the whole plane or of the perfect slip boundary condition, where there is no vorticity production near the boundary and 
maximum principle is directly applied to derive an a priori bound of the vorticity field. Although the analysis of the vorticity equations is a core part also in the proof of our Liouville theorem, 
the key idea to overcome the difficulty is to use the {\it boundary condition} on the vorticity field, rather than maximum principle.

Roughly speaking, our Liouville theorem requires four kinds of assumptions. 
The first one is a uniform bound on the velocity field including their derivatives. 
The second one is on a structure of the pressure field, which is essential to exclude the shear flows in \cite{SeSv2, Gi2} 
but is a natural requirement in order to restrict our solutions to mild solutions, i.e., solutions to the integral equations associated with  \eqref{eq.NS.intro} - \eqref{eq.noslip.intro}.  
The third one is the type I temporal decay of the velocity field as $t\rightarrow -\infty$. 
The last one is the nonnegativity of the vorticity field.  
Precisely, the main result of this paper is stated as follows.

\begin{thm}\label{thm.intro.1} Let $(u,p)$ be a solution to \eqref{eq.NS.intro}-\eqref{eq.noslip.intro} satisfying the following conditions.

\vspace{0.3cm}

\noindent {\rm (C1)} $\displaystyle \sup_{-\infty <t<0} \big ( \| u (t) \|_{C^{2+\mu}} + \| \partial_t u (t) \|_{C^\mu} \big ) <\infty$ ~~~ for some $\mu\in (0,1)$.

\noindent {\rm (C2)} ~$p=p_F + p_H$, where $p_F(t)$ is the solution to  \eqref{eq.poisson} in Proposition \ref{prop.poisson} with $F= - u (t) \otimes u (t)$
and $p_H (t)$ is the solution to \eqref{eq.laplace} in Proposition \ref{prop.laplace} with $g=\omega (t)|_{x_2=0}$, respectively.

\noindent {\rm (C3)} $\displaystyle \sup_{-\infty<t<0} (-t)^{1/2} \| u(t) \|_\infty <\infty$.

\noindent {\rm (C4)} ~$\omega\geq 0$\, in   $(-\infty,0)\times \R^2_+$, where $\omega=\partial_1u_2-\partial_2u_1$ is the vorticity field.

\vspace{0.3cm}

\noindent Then $u$ is identically zero.

\end{thm}

\noindent Here $\|\cdot\|_{C^{2+\mu}}$ and $\|\cdot\|_{C^\mu}$ denote the norms of the H${\rm \ddot{o}}$lder spaces (the definitions are stated in the end of this section), 
and $\|\cdot\|_\infty$ stands for the usual sup norm in the $x$ variables.

The condition (C3) in Theorem \ref{thm.intro.1} is compatible with the type I blow-up  assumption for forward solutions.
The sign condition (C4) on the vorticity field is a rather strong requirement at least in the class of spatially decaying solutions. 
Indeed, if there is a time $t$ such that $\sup_{x_1} |u_1(t,x_1,x_2)|\rightarrow 0$ as $x_2\rightarrow \infty$ then it is not difficult 
to see $u=0$ even when (C3) is absent; see \cite[Theorem 3.3]{Gi2}.
However, in the framework of nondecaying solutions the situation is different and becomes complicated. 
We note that, as is observed in \cite{Gi2}, there is a shear flow satisfying all of (C1), (C3), and (C4).

The key idea of the proof of Theorem \ref{thm.intro.1} is to focus on the velocity field  {\it formally} defined by the Biot-Savart law:
\begin{equation}
v (t,x) := \frac{1}{2\pi} \int_{\R^2_+} \big ( \frac{ (x-y)^\bot}{|x-y|^2} -\frac{(x-y^*)^\bot}{|x-y^*|^2} \big ) \omega (t,y) \dd y,~~~x^\bot = (-x_2,x_1), ~y^*=(y_1,-y_2).\label{eq.bs.law.intro}
\end{equation}
We note that $v$ coincides with $u$ when $u$ and $\omega$ decay fast enough at spatial infinity. By formally taking the boundary trace of $v_1$ we observe that 
\begin{equation}
v_1 (t,x_1,0) =  \frac{1}{\pi} \int_{\R^2+} \frac{y_2}{(x_1-y_1)^2 + y_2^2} \omega (t,y) \dd y. \label{trace.v_1.intro}
\end{equation}
Hence, if $v_1$ satisfies the no-slip boundary condition then the assumption (C4) implies $\omega=0$, which leads to $u=0$ by the Liouville theorem for bounded harmonic functions.

In order to justify the above formal argument we need to prove the following two claims:

\noindent Claim 1: The integral representation of the right-hand side of \eqref{eq.bs.law.intro} is well-defined. In other words, the vorticity field has an enough spatial decay so that the integral in \eqref{eq.bs.law.intro} converges.

\noindent Claim 2: The tangential component $v_1$ satisfies the no-slip boundary condition. 

Both of two claims are far from trivial, for we have to start from the spatially nondecaying data, and the right-hand side of  \eqref{trace.v_1.intro} is highly nonlocal. To show Claim 1 we make use of the type I temporal decay of $u$ assumed in (C3). In fact, since (C3) is a scaling invariant bound, by applying the result of \cite{CarLo} or \cite{MaTo} we can establish the Gaussian pointwise bound of the Green function for the heat-transport operator $\partial_t - \Delta + u\cdot \nabla$ with the Neumann boundary condition. This pointwise estimate of the Green function leads to a polynomial decay of the vorticity field as $x_2\rightarrow \infty$, which makes the integral of \eqref{eq.bs.law.intro} well-defined. The key ingredient of the proof of Claim 2 is the boundary condition on the vorticity field. Indeed,  combined with a calculation based on the integration by parts, the vorticity boundary condition yields $\partial_t v_1 (t,x_1,0)=0$ for $-\infty<t<0$ and $x_1\in \R$, as is already observed in \cite{Mae} in the setting of spatially decaying solutions. Then the no-slip boundary condition for $v_1$ is a consequence of the convergence $\lim_{t\rightarrow -\infty} v_1 (t,x_1,0)=0$, which can be verified from the time decay condition (C3) and the polynomial decay of the vorticity field established in Claim 1. 
 
As an application of Theorem \ref{thm.intro.1}, we can extend the geometric regularity criterion in \cite {GiMi} for the three-dimensional Navier-Stokes equations in the half space to the case of the no-slip boundary condition.

\begin{thm}\label{thm.intro.2}
Let $(u,p)$  be a spatially bounded mild solution to the Navier-Stokes equations \eqref {eq.NS.Sec.6}-\eqref {eq.Dirichlet.Sec.6} in $(0, T) \times \R^3_+$. Assume that the possible blow-up of $u$ is type I, i.e.
\begin{align*}
\sup_{0<t<T}   (T-t)^\frac12 \|u(t) \|_\infty <\infty.
\end{align*}
Let $d$ be a positive number and let $\eta$ be a nondecreasing continuous function on $[0,\infty)$ satisfying $\eta(0)=0$. 
Assume that $\eta$ is a modulus of continuity in the $x$ variables for the vorticity direction $\xi=\omega/|\omega|$, in the sense that 
\begin{align}\tag{CA}
|\xi (t,x)-\xi (t,y)|\leq \eta (|x-y|)~~~~ {\rm for} ~~ (t,x),(t,y)\in {\Omega}_d,
\end{align}
where $\Omega_d = \{(t,x) \in (0,T)\times \R^3_+ \mid |\omega (x,t)|>d\}$.
Then $u$ is bounded up to $t=T$.
\end{thm}
 
The condition {\rm (CA)} is called a ``continuous alignment'' condition.
This kind of geometric condition on the vorticity direction was firstly given in \cite{CoFe} for a finite energy solution in $\R^3$ with $H^1$ initial data.
In \cite{CoFe} the modulus $\eta$ is taken as $\eta(\sigma) = A\sigma$ with some constant $A>0$, while the type I condition is not needed there. 
The condition in \cite{CoFe} was relaxed in \cite{BaBe}, where $\eta$ is allowed to be $\eta(\sigma)=A\sigma^{1/2}$; see \cite{GiMi} for further references on the related results. 
A corresponding result to \cite{BaBe} for slip boundary conditions is established in \cite{B2}, where $\eta (\sigma) = A \sigma^{1/2}$ in (CA).
However, under the no-slip boundary condition the regularity criterion, so far obtained in \cite{B3}, needs an extra assumption that the boundary integral of
the normal derivative of the square of the vorticity is sufficiently small.
As far as the authors know, the present paper gives the first contribution to the case of the no-slip boundary condition
under the same assumption to the whole space. 
This is rather surprising since the geometric regularity criterion is still valid even if the vorticity is created from the boundary
because of the no-slip boundary condition.  As in \cite{GiMi}, the proof of Theorem \ref{thm.intro.2} is based on a blow-up argument.

Before concluding this section, we introduce Banach spaces with nondecaying functions. Let $\Omega$ be a domain in $\R^n$, $n\in \N$. Then, for $k\in \N\cup \{0\}$ and $\mu\in (0,1)$ the spaces $BC(\overline{\Omega})$, $C^k(\overline{\Omega})$, and $C^{k+\mu}(\overline{\Omega})$  are respectively defined by 
\begin{align*}
& BC (\overline{\Omega}) = \big \{\, f \in C(\overline{\Omega})~|~\| f\|_\infty = \sup_{x\in \overline{\Omega}} |f(x)| <\infty \, \big\},\\
& C^k (\overline{\Omega}) = \big \{ f\in BC (\overline{\Omega})~|~ \nabla^\alpha f\in BC(\overline{\Omega}), \, |\alpha|\leq k,~~~ \| f \|_{C^k} =  \sum_{|\alpha|\leq k} \| \nabla^\alpha f \|_\infty<\infty \, \big \},\\
& C^{k+\mu} (\overline{\Omega}) = \big  \{ f\in C^k (\overline{\Omega})~|~\\
& ~~~~~~~~~~~~~~~ \| f\|_{C^{k+\mu}} = \|  f \|_{C^k} + \sum_{|\alpha|=k} \sup_{x,y\in \overline{\Omega},\,x\ne y}\frac{|\nabla^\alpha f(x)-\nabla^\alpha f(y)|}{|x-y|^\mu} <\infty \,  \big \}.
\end{align*}
Let us also introduce the $BMO$ spaces as follows.
\begin{align*}
& BMO (\R^n) = \big \{ f\in L^1_{loc}(\R^n)~|~\|f \|_{BMO} = \sup_{B} \frac{1}{|B|} \int_B |f - {\rm Avg}_B f | \dd x <\infty \,\big \},\\
& BMO(\Omega) = \big \{ f\in L^1_{loc}(\Omega) ~|~ {\rm there~is~} g\in BMO(\R^n)~{\rm such~that~} f =g ~~{\rm a.e.~in}~\Omega,\\
& ~~~~~~~~~~~~~~~~ \|f \|_{BMO}=\inf \{\|g\|_{BMO}~|~g\in BMO(\R^n),~f=g~~{\rm a.e.~in ~}\Omega \} \, \big\}. 
\end{align*}
In the definition of $\|\cdot\|_{BMO}$ the supremum is taken over all ball $B$ in $\R^n$, $|B|$ is the volume of $B$, and ${\rm Avg}_B  f = |B|^{-1}\int_B f \dd x$. 

This paper is organized as follows. In Section \ref{sec.bc.stokes} we consider the Stokes equations with a inhomogeneous term and derive the boundary condition on the vorticity field. We also obtain the integral equations for the vorticity field, which is useful to estimate the vorticity field directly. Section \ref{sec.liouville} is the core part of  this paper, and we study \eqref{eq.NS.intro} - \eqref{eq.noslip.intro} under the conditions of Theorem \ref{thm.intro.1}. To this end we establish a temporal decay estimate in Section \ref{sec.temporal.decay} and a spatial decay estimate in Section \ref{sec.spatial.decay}. Claim 1 and Claim 2 in this section are respectively stated as Lemma \ref{lem.spatial.decay.vorticity} and Lemma \ref{lem.biot-savart}. These are proved in Section \ref{sec.bs.law},  which completes the proof of Theorem \ref{thm.intro.1}.  Finally,  as an application of Theorem \ref{thm.intro.1}, we prove Theorem \ref{thm.intro.2} in Section \ref{sec.criterion}.

After this work was completed, the result was presented by the first author
in the Clay workshop ``The Navier-Stokes equations and related topics'' on
September 29, 2013. The authors are grateful to Professor Gregory Seregin who
kindly pointed out during the workshop that a Liouville type result can
be proved without using the vorticity equation \cite{Sepre}. However,
his results need an assumption that the (kinetic) energy is bounded in
time. This assumption imposes a decay at the spatial infinity and it is
not enough to apply for proving a geometric regularity criterion such as Theorem \ref{thm.intro.2}.

The work of the first author is partly supported by the Japan Society of
the Promotion of Science (JSPS) through grants Kiban (S) 21224001, Kiban
(A) 23244015 and Houga 25610025. The work of the third author is partly
supported by JSPS through Grant-in-Aid for Young Scientists (B) 25800079.

\section{Vorticity boundary condition for the Stokes flows}\label{sec.bc.stokes}

In this section we consider the Stokes equations
\begin{align}
\partial_t u - \Delta u + \nabla p = {\rm div}\, F,~~~~~ {\rm div}\, u =0~~~~~~~~~{\rm in}~~ (-L, 0) \times \R^2_+\label{eq.Stokes}
\end{align}
subject to the no-slip boundary condition
\begin{align}
u=0~~~~~~~~~~{\rm on} ~~ (-L, 0) \times \partial \R^2_+.\label{eq.noslip}
\end{align}

The aim of this section is to derive the boundary condition on the vorticity field
\[
\omega = -\nabla^\bot \cdot u,~~~~~~\nabla^\bot = (\partial_2,-\partial_1)^\top.
\]
If the flow possesses enough spatial decay then the vorticity boundary condition can be derived  from the Biot-Savart law (e.g. \cite{Mae}). Here we give an alternative derivation of the vorticity boundary condition in order to deal with  nondecaying  flows. The derivation is closely related with the structure of the pressure field.  As is well-known, by acting the {\rm div} operator in \eqref{eq.Stokes}  the pressure field is recovered as a solution to the Poisson equations with the inhomogeneous Neumann boundary condition. With this in mind we start from

\begin{prop}\label{prop.poisson} Assume that $F= (F_{ij})_{1\leq i,j\leq 2} \in (C^2 (\overline{\R^2_+}))^{2\times 2}$, $F_{ij}=\partial_2 F_{ij}=0$ on $\partial\R^2_+$ for each $i,j$. Then there is a unique (up to a constant) solution $p_F \in BMO (\R^2_+)$ to
\begin{equation}\label{eq.poisson}
\begin{cases}
& \Delta p_F = {\rm div} \, {\rm div}\, F ~~~~~~~{\rm in}~~ \R^2_+,\\
&  \partial_2 p_F =0~~~~~~~~~~ ~~~~~~{\rm on}~~ \partial\R^2_+,
\end{cases}
\end{equation}
such that
\begin{align}
& \| p_F\|_{BMO} \leq C \| F \|_{\infty},~~~~~~~~\| \nabla p_F \|_{C^\mu} \leq C \| F \|_{C^{1+\mu}},~~~~~0<\mu < 1.\label{est.prop.poisson}
\end{align}

\end{prop}

\noindent {\it Proof.} As usual, let us introduce  the even extension: $\tilde p_F (x) = p_F (x) $ for $x_2\geq 0$ and $\tilde p_F (x) = p_F (x^*)$ for $x_2<0$. The same extension is introduced also for $F_{11}$ and $F_{12}$, while the odd extension is applied for $F_{12}$ and $F_{21}$.  We denote by $\tilde F$ the tensor extended in this manner.  Then \eqref{eq.poisson} is reduced to the Poisson equation $\Delta \tilde p_F = {\rm div}\,{\rm div}\, \tilde F$ in $\R^2$ by the assumption $F_{ij}=\partial_2 F_{ij}=0$ on $\partial \R^2_+$. Its  solution is written as $\tilde p_F = - {\rm div}\, {\rm div}\, (-\Delta_{\R^2})^{-1} \tilde F$, where the operator $(-\Delta_{\R^2})^{-1}$ is defined as the convolution with the Newton potential in $\R^2$. It is well known that ${\rm div}\, {\rm div}\, (-\Delta_{\R^2})^{-1}$ defines a singular integral operator, and hence it is bounded in $BMO (\R^2)$, and $\nabla{\rm div}\, {\rm div}\, (-\Delta_{\R^2})^{-1}$ is bounded from $C^{1+\mu} (\R^2)$ to $C^\mu (\R^2)$, $0<\mu<1$. Thus  \eqref{est.prop.poisson} holds. The uniqueness is a consequence of the classical Liouville theorem for harmonic functions in $\R^2$. The proof is complete.

\vspace{0.5cm}

In order to recover the no-slip boundary condition on the velocity field we need to introduce the harmonic pressure field. As a preliminary, let us recall some results on the fractional power of the Laplace operator $-\partial_1^2$. As is well known, $-\partial_1^2$ is realized as a sectorial operator in $BC(\R)$ (cf. \cite{Lu}), and hence its fractional power $(-\partial_1^2)^{1/2}$ is also sectorial in $BC(\R)$. The characterization of the interpolation spaces as in \cite[Section 3.1.3]{Lu} implies that
\begin{equation}
C^{1+\mu} (\R) \hookrightarrow D( (-\partial_1^2)^\frac12 ) ~~~~~~~~{\rm for~all}~~\mu \in (0,1),\label{embedding.1/2}
\end{equation}
where $D( (-\partial_1^2)^{1/2} )$ is the domain of $(-\partial_1^2)^{1/2}$ in $BC(\R)$. Note that the semigroup generated by $(-\partial_1^2)^{1/2}$ is nothing but the Poisson semigroup whose kernel is explicitly described.

\begin{prop}\label{prop.laplace} Assume that $g\in C^{1+\mu} (\R)$ for some $\mu\in (0,1)$. Then there is a unique (up to a constant) solution $p_H\in  L^1_{loc} (\R^2_+)$  to
\begin{equation}\label{eq.laplace}
\begin{cases}
& \Delta p_H = 0~~~~~~~~~~~~{\rm in}~~~ \R^2_+,\\
&  \partial_2 p_H = \partial_1 g~~~~~~~~{\rm on}~~ \partial\R^2_+,
\end{cases}
\end{equation}
such that
\begin{align}
\sup_{x\in \R^2_+} x_2 | \nabla p_H (x) |  \leq C \| g \|_{\infty},~~~~~~~\| \nabla p_H \|_{C^{\mu'}} \leq C \| g \|_{C^{1+\mu}},~~~0<\mu'<\mu.\label{est.prop.laplace.1}
\end{align}
Moreover, it follows that
\begin{equation}
\lim_{x_2\downarrow 0} \partial_1 p_H (x) = (-\partial_1^2 )^\frac12 g (x_1) ~~~~~~~{\rm in} ~~BC (\R).\label{est.prop.laplace.2}
\end{equation}
\end{prop}

\begin{rem}{\rm  In Proposition \ref{prop.laplace} the weight estimate in \eqref{est.prop.laplace.1} is essential in view of the uniqueness of solutions. In particular,  if one tries  to avoid the Poiseuille type flows as in \cite{Gi2}  it is important to impose suitable conditions on the behavior of the harmonic pressure at spatial infinity.

}
\end{rem}

\noindent {\it Proof of Proposition \ref{prop.laplace}.} The solution $p_H$ is constructed so as to satisfy the representation
\begin{align}
\nabla p_H (x) = - \int_0^\infty \big ( \nabla \partial_1 e^{- (x_2 + y_2 ) (-\partial_1^2)^\frac12 }  g \big ) (x_1) \dd y_2. \label{proof.prop.laplace.1}
\end{align}
Indeed, if $g$ is compactly supported then $p_H$ is given by $p_H= - \int_0^\infty \partial_1 e^{-(x_2+y_2) (-\partial_1^2)^{1/2} } g \dd y_2$, where the integral converges absolutely. Then we modify $p_H$ by adding a constant so that  the condition $p_H (0)=0$  holds and both \eqref{eq.laplace} and  \eqref{proof.prop.laplace.1} are satisfied. We denote this modified solution by $p_H (g)$. The straightforward calculation of the Poisson semigroup yields that
\begin{align}
\| \nabla^k e^{- x_2 (-\partial_1^2)^\frac12 }  g \|_{L^\infty_{x_1}} \leq C x_2^{-k} \| g\|_{\infty},~~~~~~~k=0,1,2,\label{proof.prop.laplace.2}
\end{align}
and
\begin{align}
\| \nabla p_H (g) \|_{C^{\mu'}} \leq C \| g \|_{C^{1+\mu}},~~~~~0<\mu'<\mu<1.\label{proof.prop.laplace.3}
\end{align}
Then for general $g\in C^{1+\mu}(\R)$ we approximate $g$ by $g\chi_R$ with a smooth cut-off $\chi_R$ and take the limit of $p_H (g_R)$ at $R\rightarrow \infty$. Since $g_R \rightarrow g$ in $C^1 (K)$ for each compact set $K\subset \R$ and $\sup_{R>0}\| g_R \|_{C^{1+\mu}}<\infty$, it is not difficult to show that there is a  subsequence of $\{p_H (g_R)\}_{R>0}$ which converges to some $p_H$ in $C^1 (K')$ for each compact set $K'\subset  \overline{\R^2_+}$. It is easy to see that $p_H$ solves \eqref{eq.laplace} and also satisfies \eqref{proof.prop.laplace.1} by the Lebesgue convergence theorem.  The estimate \eqref{est.prop.laplace.1} is a consequence of \eqref{proof.prop.laplace.2} and \eqref{proof.prop.laplace.3}. To show \eqref{est.prop.laplace.2} we observe from \eqref{proof.prop.laplace.1} that
\begin{align*}
\partial_1 p_H (x)  = \int_0^\infty (-\partial_1^2) e^{-(x_2 + y_2 ) (-\partial_1^2)^\frac12} g \dd y_2 & =  - \int_0^\infty  (-\partial_1^2)^\frac12 \partial_{y_2} \big ( e^{-(x_2 + y_2 ) (-\partial_1^2)^\frac12} g \big ) \dd y_2\\
& = (-\partial_1^2)^\frac12 e^{-x_2  (-\partial_1^2)^\frac12} g,~~~~~~{\rm for}~~x_2>0.
\end{align*}
Hence \eqref{est.prop.laplace.2} follows from \eqref{embedding.1/2}. The uniqueness of solutions to \eqref{eq.laplace} is again reduced to the classical Liouville theorem for harmonic functions in $\R^2$ by a suitable reflection argument. The details are omitted here. The proof is now complete.

\vspace{0.5cm}

We are now in position to derive the vorticity boundary condition for nondecaying flows.

\begin{lem}\label{lem.boundary.vorticity} Assume that $F= (F_{ij})_{1\leq i,j\leq 2} \in C((-L,0)\times (C^2 (\overline{\R^2_+}))^{2\times 2} )$, $F_{ij} (t) =\partial_2 F_{ij} (t) =0$ on $(-L,0)\times \partial\R^2_+$ for each $i,j$.  Let $(u,p)$ be the solution to \eqref{eq.Stokes}-\eqref{eq.noslip} such that

\vspace{0.3cm}

\noindent {\rm (C1)} $\displaystyle \sup_{-L<t<0} \big ( \| u (t) \|_{C^{2+\mu}} + \| \partial_t u (t) \|_{C^\mu}  \big ) <\infty$ ~~~ for some ~$\mu\in (0,1)$,

\noindent {\rm (C2)} $p=p_F + p_H$, where $p_F(t)$ is the solution to  \eqref{eq.poisson} in Proposition \ref{prop.poisson} with $F=F(t)$ and $p_H (t)$ is the solution to \eqref{eq.laplace} in Proposition \ref{prop.laplace} with $g=\omega (t)|_{x_2=0}$, respectively.

\vspace{0.3cm}

\noindent Then $\omega$ satisfies
\begin{align}
\partial_t \omega  - \Delta \omega = -\nabla^\bot \cdot {\rm div} \, F ~~~~~~~~~~ {\rm in}~~ (-L, 0) \times \R^2_+\label{eq.vorticity.linear}
\end{align}
with
\begin{align}
\partial_2 \omega     + (-\partial_1^2 )^\frac12 \omega = - \partial_1 p_F~~~~~~~~ {\rm on}~~ (-L, 0) \times \partial \R^2_+\label{eq.vorticity.boundary.linear}.
\end{align}

\end{lem}

\noindent {\it Proof.} It is straightforward to see \eqref{eq.vorticity.linear}. To show \eqref{eq.vorticity.boundary.linear} we first recall the equality $-\Delta u = \nabla^\bot \omega$ and then  \eqref{eq.Stokes} yields $\partial_2 \omega = -\partial_t u_1 -\partial_1 p + {\bf \tau} \cdot {\rm div}\, F$ for $x_2>0$, where ${\bf \tau} = (1,0)^\top$. Thus we have
\begin{align*}
\lim_{x_2 \downarrow 0} \partial_2 \omega & = -  \lim_{x_2 \downarrow 0} \partial_t u_1 - \lim_{x_2 \downarrow 0} \partial_1 p_F - \lim_{x_2 \downarrow 0}\partial_1 p_H +  \lim_{x_2 \downarrow 0}{\bf \tau} \cdot {\rm div}\, F\\
& = -\partial_1 p_F|_{x_2=0} - (-\partial_1^2)^\frac12 \omega |_{x_2=0}~~~~~~~~~{\rm by}~~ \eqref{est.prop.laplace.2}.
\end{align*}
The proof is now complete.

\vspace{0.5cm}

Lemma \ref{lem.boundary.vorticity} leads to the integral equation for the vorticity field, which is useful to estimate the vorticity directly including near the boundary.

Let $G(t,x) = (4\pi t ) ^{-1} \exp \big ( -|x|^2/ (4 t) \big )$ be the two-dimensional Gaussian. Then for each $t>0$ we introduce the operator $e^{tB}$ defined by
\begin{equation}
e^{t B} f = G(t) * f + G (t) \star f  + \Gamma (t)\star f,\label{def.e^B}
\end{equation}
where
\begin{equation}
\Gamma (t) = 2 \int_0^\infty  \big (\partial_1^2 + (-\partial_1^2 )^\frac{1}{2} \partial_2 \big )  G (t+\tau ) \dd\tau \label{def.Gamma}
\end{equation}
with the notations
\begin{equation*}
f ~* ~ h (x) = \int_{\R^2_+} f (x-y) h(y) \dd y, ~~~~~~ f ~\star ~ h (x ) = \int_{\R^2_+} f (x-y^*) h(y) \dd y, ~~~~~y^* = (y_1,-y_2 ).\label{def.star}
\end{equation*}
For $g\in C_0^\infty (\R)$ we set  
\[
e^{tB} (g\delta_{\partial\R^2_+}) = \int_\R K(t,x,y) |_{y_2=0}  g(y_1) \dd y_1, 
\]
where $K(t,x,y)$ is the kernel of $e^{tB}$. Due to the pointwise estimate of $K(t,x,y)$ in \eqref{est.prop.pointwise.kernel}, the term $e^{tB} (g\delta_{\partial\R^2_+})$ makes sense also for $g\in L^\infty (\R)$. The operator $e^{tB}$ naturally arises in the vorticity equations. Indeed, if $f\in C_{0}^\infty (\R^2_+)$ then $e^{tB} f$ satisfies the (homogeneous) vorticity equations \eqref{eq.vorticity.linear}-\eqref{eq.vorticity.boundary.linear}, i.e., $w (t) = e^{(t+L)B} f$ solves
\begin{equation}
\partial_t w - \Delta w =0 ~~~~ {\rm in}~~ (-L,0)\times \R^2_+, ~~~~~~~ \partial_2 w + (-\partial_1^2)^\frac12 w =0 ~~~~ {\rm on}~~ (-L,0) \times \partial\R^2_+,
\end{equation}
but with the initial data $w (-L)= \displaystyle \lim_{t\rightarrow -L} w (t) =f + \Gamma (0)\star f$ in $L^p(\R^2_+)$ for all $1<p<\infty$; see \cite[Sections 3,5]{Mae} for details. Note that $\Gamma (0) \star$ is a singular integral operator. In particular, we have $\|\Gamma (0)\star f \|_{L^p}\leq C \| f \|_{L^p}$ for all $f\in L^p (\R^2_+)$. If $f=-\nabla^\bot \cdot u$ with $u\in C_{0,\sigma}^\infty (\R^2_+)$ then $\Gamma (0)\star f=0$ (see \cite[Proposition 3.2]{Mae}), hence in this case we recover the initial condition $w (-L) =f$, as desired. For each $t>0$ let us introduce the operator $T(t): (L^\infty (\R^2_+))^2 \rightarrow L^\infty (\R^2_+)$ as follows:
\begin{align}
\langle T(t) v, f\rangle _{L^2} = \langle v_1, \partial_2 e^{t B} f \rangle _{L^2}   - \langle v_2, \partial_1 e^{t B} f\rangle _{L^2}~~~~~~{\rm for~all}~~f\in L^1 (\R^2_+). \label{def.T(t)}
\end{align}
Here $\langle ,\rangle _{L^2}$ denote the inner product of $L^2 (\R^2_+)$. The operator $T(t)$ is well-defined due to the estimate $\| \nabla e^{t B} f \|_{L^1} \leq C t^{-1/2} \| f\|_{L^1}$ by \cite[Lemma 3.4]{Mae} and the duality $L^1 (\R^2_+)^*=L^\infty (\R^2_+)$. In particular, we have
\begin{align}
\| T(t) v \|_{\infty} \leq C t^{-\frac12} \| v\|_{\infty},~~~~~~~~t>0.\label{est.T.1}
\end{align}

\begin{lem}\label{lem.integral} Assume that the conditions in Lemma \ref{lem.boundary.vorticity} hold and ${\rm div}\, F \in ( L^\infty (\R^2_+))^2$.  Then $\omega$ satisfies the integral equation
\begin{align}
\omega (t) = T (t-s) u (s) + \int_s^t T(t-\tau ) {\rm div}\, F (\tau ) \dd \tau + \int_s^t e^{(t-\tau)B} (\partial_1 p_{F(\tau )}\delta_{\partial\R^2_+})\dd \tau \label{eq.lem.integral}
\end{align}
for $-L<s<t<0$.
\end{lem}

\begin{rem}{\rm There are several solution formulas for the velocity field of the Stokes flows
in the half space with the Dirichlet condition for example in \cite{So,U}.
Ours differs from those in the literature since it is a convenient form to represent the vorticity field.

}
\end{rem}

\noindent {\it Proof of Lemma \ref{lem.integral}.} Take any $\phi (\tau ,x)\in C_{0}^\infty ([s,t]\times
\overline{\R^2_+})$. Multiplying \eqref{eq.vorticity.linear} by $\phi$
 and using the integration  by parts, we observe that $\omega$ satisfies
\begin{align*}
\langle \omega (t), \phi (t)\rangle _{L^2} &=\int_{\R^2_+} u(s)\cdot \nabla^\bot \phi (s) \dd x + \int_s^t \int_{\R^2_+} {\rm div}\, F \cdot \nabla^\bot \phi (\tau) \dd  x \dd \tau \\
&~~~- \int_s^t \int_{\partial \R^2_+} (\phi \partial_2 \omega - \omega \partial_2 \phi)(\tau) \dd x_1 \dd \tau + \int_s^t \int_{\R^2_+}  \omega (\partial_\tau \phi +  \Delta \phi ) (\tau) \dd x \dd \tau.
\end{align*}
Fix $R\gg 1$ and set $\phi_R (\tau ,x):= ( \chi_R e^{(t-\tau)B}\psi ) (x)$, where $\psi \in C_{0}^\infty (\R^2_+)$ and $\chi_R=\chi_R(x)$ is a nonnegative smooth cut-off function in $\R^2$ such that $\chi_R(x)=1$ if $|x|\leq R$ and $\chi_R(x)=0$ if $|x|\geq 2 R$. We may assume that $\|\nabla^k \chi_R\|_\infty\leq C R^{-k}$ for $k=0,1,2$. Then we set
\begin{align*}
\langle \omega (t), \phi_R (t)\rangle _{L^2} &=\int_{\R^2_+} u(s)\cdot \nabla^\bot \phi_R (s) \dd x + \int_s^t \int_{\R^2_+} {\rm div}\, F  \cdot \nabla^\bot \phi_R (\tau) \dd  x \dd \tau \\
&~- \int_s^t \int_{\partial \R^2_+} (\phi_R \partial_2 \omega - \omega \partial_2 \phi_R)(\tau) \dd x_1 \dd \tau + \int_s^t \int_{\R^2_+} \omega (\partial_\tau \phi_R + \Delta \phi_R)  (\tau) \dd x \dd \tau\\
&:= I_1 + I_2 - I_3 + I_4.
\end{align*}
As for $I_1$, we have 
\begin{align}
I_1 &=\int_{\R^2_+} u(s)\cdot \nabla^\bot (\chi_R e^{(t-s)B}\psi ) \dd x \nonumber \\
&=\int_{\R^2_+} \big ( u_1(s)  (\partial_{2}\chi_R) e^{(t-s)B}\psi-u_2(s) (\partial_{1}\chi_R) e^{(t-s)B}\psi  \big )\dd x \nonumber \\ 
&~~~ +\int_{\R^2_+}\chi_R \big ( u_1(s)\partial_{2}e^{(t-s)B}\psi - u_2(s)\partial_{1}e^{(t-s)B}\psi \big )\dd x.\label{eq.I_1}
\end{align}
Thanks to \cite[Lemma 3.4]{Mae} and (C1) we have  $e^{(t-s)B}\psi \in L^p(\R^2_+)$ for any $1<p \leq \infty$, and the first term of right-hand side of \eqref{eq.I_1} converges to zero in the limit  $R\rightarrow \infty$. As for the second term of \eqref{eq.I_1}, we observe from \cite[Lemma 3.4]{Mae} that $\|\nabla e^{(t-s)B} \psi\|_{L^1}\leq C (t-s)^{-1/2} \|\psi \|_{L^1}$. Hence the H${\rm \ddot{o}}$lder  inequality implies that $u_1(s)\partial_{2}e^{(t-s)B}\psi - u_2(s)\partial_{1}e^{(t-s)B}\psi$ belongs to $L^1 (\R^2_+)$ for $t>s$. Thus we have $\lim_{R\rightarrow \infty}I_1=\langle T (t-s) u (s), \psi  \rangle _{L^2}$ by the definition of $T(t-s)$. Similarly, by the assumption ${\rm div}\, F\in ( L^\infty (\R^2_+))^2$ and by the Fubini theorem, we have $\lim_{R\rightarrow \infty}I_2=\langle \int_s^t T(t-\tau ) {\rm div} \, F (\tau ) \dd \tau, \psi \rangle _{L^2}$.

As for $I_3$, we recall the vorticity boundary condition \eqref{eq.vorticity.boundary.linear}. Then it follows that
\begin{align*}
I_3&= \int_s^t \int_{\partial \R^2_+} (\phi_R \partial_2 \omega - \omega \partial_2 \phi_R)(\tau) \dd x_1 \dd \tau\\
&=\int_s^t \int_{\partial \R^2_+} \big \{ - \chi_R \big (\partial_1 p_{F(\tau )}+ (-\partial_1^2 )^\frac12 \omega \big ) e^{(t-\tau)B}\psi \\
& ~~~~~~~~~~~~~~ - \omega \big ( \chi_R \partial_{2}e^{(t-\tau)B}\psi +(\partial_{2}\chi_R)e^{(t-\tau)B}\psi  \big ) \big \}\dd x_1 \dd \tau\\
&= - \int_s^t \int_{\partial \R^2_+}  \chi_R \big (  \partial_1 p_{F(\tau )}  e^{(t-\tau)B}\psi  + \omega   (-\partial_1^2 )^\frac12e^{(t-\tau)B}\psi  + \omega \partial_{2}e^{(t-\tau)B}\psi \big ) \dd x_1 \dd \tau \\
& ~~~ - \int_s^t \int_{\partial \R^2_+}  (\partial_{2}\chi_R)  \omega e^{(t-\tau)B}\psi  \dd x_1 \dd \tau.
\end{align*}
Since $ (-\partial_1^2 )^\frac12 e^{(t-\tau)B}\psi + \partial_{2}e^{(t-\tau)B}\psi =0$ on $\partial \R^2_+$, we obtain
\begin{align*}
\lim_{R\rightarrow \infty}I_3=- \lim_{R\rightarrow \infty}\int_s^t \int_{\partial \R^2_+}\chi_R \partial_1 p_{F(\tau )} e^{(t-\tau)B}\psi  \dd x_1 \dd \tau =- \langle \int_s^t e^{(t-\tau)B}\partial_1 p_{F(\tau )}\delta_{\partial\R^2_+}, \psi  \rangle _{L^2}.
\end{align*}
Finally, we consider  $I_4$. It is easy to check that
\begin{align*}
\partial_{\tau}\phi_R + \Delta \phi_R = (\Delta \chi_R ) e^{(t-\tau)B}\psi +2 \nabla \chi_R \cdot \nabla e^{(t-\tau)B}\psi.
\end{align*}
Since $\omega$ is bounded in space and time, by using $\|\nabla^k \chi_R\|_\infty\leq C R^{-k}$ for $k=1,2$ and the estimate of $e^{(t-\tau)B} \psi$ the term $I_4$ is shown to converge to zero as $R\rightarrow \infty$.
Combining the above calculations, we have
\begin{align*}
&\lim_{R\rightarrow \infty}\langle \omega (t), \phi_R (t)\rangle _{L^2}=\lim_{R\rightarrow \infty}[I_1+I_2-I_3+I_4]\\
&=\langle T (t-s) u (s) + \int_s^t T(t-\tau ) {\rm div} \, F (\tau ) \dd \tau + \int_s^t e^{(t-\tau ) B} \big ( \partial_1 p_{F(\tau )} \delta_{\partial\R^2_+} \big ) \dd \tau, \psi  \rangle _{L^2}
\end{align*}
Note that $\phi_R (t)=\lim_{\tau \rightarrow t} \chi_R e^{(t-\tau)B}\psi =\chi_R (\psi  + \Gamma (0)\star \psi )$. By the definition of $\Gamma (0)$ we have $(\partial_2 + (-\partial_1^2)^{1/2} ) \Gamma (0)\star \psi=0$ in $\R^2_+$. Then, together with the divergence free property of $u$ and $u=0$ on $\partial\R^2_+$, we observe from the integration by parts  and $\| \Gamma (0)\star \psi \|_{L^p} + \| \partial_1 (-\partial_1^2)^{-\frac12} \Gamma (0)\star \psi \|_{L^p} \leq C \| \psi \|_{L^p}$ that 
\begin{align*}
\lim_{R\rightarrow \infty} \langle u_1 , \chi_R \partial_{2} \Gamma (0)\star \psi \rangle_{L^2}  & =  - \lim_{R\rightarrow \infty} \langle u_1 , \chi_R  (-\partial_1^2)^{\frac12} \Gamma (0)\star \psi \rangle_{L^2}  \\
&=  - \lim_{R\rightarrow \infty} \langle \partial_1 u_1 , \chi_R  \partial_1 (-\partial_1^2)^{-\frac12} \Gamma (0)\star \psi \rangle_{L^2}  \\
& = - \lim_{R\rightarrow \infty} \langle u_2 , \chi_R  \partial_1 (-\partial_1^2)^{-\frac12} \partial_2 \Gamma (0)\star \psi \rangle_{L^2}\\
& =  \lim_{R\rightarrow \infty} \langle u_2 , \chi_R  \partial_1 \Gamma (0)\star \psi \rangle_{L^2},
\end{align*}
that is, $\lim_{R\rightarrow \infty}\langle \omega (t), \chi_R \Gamma (0)\star \psi \rangle_{L^2}=0$, again from the integration by parts for $\omega = \partial_1 u_2 - \partial_2 u_1$. Thus it follows that
\begin{align*}
\lim_{R\rightarrow \infty}\langle \omega (t), \phi_R (t)\rangle _{L^2} & =\lim_{R\rightarrow \infty}\langle \omega (t), \chi_R \psi \rangle _{L^2} + \lim_{R\rightarrow \infty}\langle \omega (t), \chi_R \Gamma (0)\star \psi \rangle_{L^2}\\
& =\lim_{R\rightarrow \infty}\langle \omega (t), \chi_R \psi \rangle _{L^2} = \langle \omega (t),  \psi \rangle _{L^2}.
\end{align*}
Since $\psi$ is arbitrary the proof is now complete.

\vspace{0.3cm}

%%%%%%%%%%%%%%%%%%%%%%%%%%%%%%%%%%%%%%%%

As an immediate consequence of Lemmas \ref{lem.boundary.vorticity} and \ref{lem.integral}, we obtain the vorticity equations for the full nonlinear problem \eqref{eq.NS.intro}-\eqref{eq.noslip.intro}.

\begin{prop}\label{prop.vorticity.equation} Let $(u,p)$ be the solution to \eqref{eq.NS.intro}-\eqref{eq.noslip.intro} such that

\vspace{0.3cm}

\noindent {\rm (C1)} $\displaystyle \sup_{-\infty<t<0} \big ( \| u (t) \|_{C^{2+\mu}} + \| \partial_t u (t) \|_{C^\mu}  \big ) <\infty$ ~~~ for some ~ $\mu\in (0,1)$,

\noindent {\rm (C2)} ~$p=p_F + p_H$, where $p_F(t)$ is the solution to  \eqref{eq.poisson} in Proposition \ref{prop.poisson} with $F= - u(t)\otimes u(t)$ and $p_H (t)$ is the solution to \eqref{eq.laplace} in Proposition \ref{prop.laplace} with $g=\omega (t)|_{x_2=0}$, respectively.

\vspace{0.3cm}

\noindent Then $\omega$ satisfies
\begin{align}
\partial_t \omega  - \Delta \omega = \nabla^\bot \cdot {\rm div}\, ( u\otimes u ) ~~~~~~~~~~~ {\rm in}~~ (-\infty, 0) \times \R^2_+\label{eq.vorticity}
\end{align}
with
\begin{align}
\partial_2 \omega     + (-\partial_1^2 )^\frac12 \omega = - \partial_1 p_F~~~~~~~ {\rm on}~~ (-\infty, 0) \times \partial \R^2_+\label{eq.vorticity.boundary}.
\end{align}
Moreover, $\omega$ satisfies the integral equation \eqref{eq.lem.integral} for $-\infty<s<t<0$.

\end{prop}

\section{Liouville type result}\label{sec.liouville}

In this section we prove Theorem \ref{thm.intro.1}. As stated in the introduction, the key idea of the proof  is to derive the spatial decay of vorticity fields in the vertical direction and to verify the relation of the Biot-Savart law  between the velocity and  the vorticity. More precisely, the core parts of the proof are  the following two lemmas.

\begin{lem}\label{lem.spatial.decay.vorticity} Under the conditions {\rm (C1)}, {\rm (C2)}, and {\rm (C3)} of Theorem \ref{thm.intro.1} the vorticity $\omega$ satisfies
\begin{equation}
\sup_{(t,x)\in (-\infty,0)\times \R^2_+} x_2^{1+\theta} |\omega (t,x) | <\infty~~~~~~~~{\rm for~all}~~ \theta\in (0,1).\label{est.lem.spatial.decay.vorticity}
\end{equation}

\end{lem}

\begin{lem}\label{lem.biot-savart} Under the conditions {\rm (C1)}, {\rm (C2)}, and {\rm (C3)} of Theorem \ref{thm.intro.1} the velocity $u$ is represented as
\begin{equation}
u (t,x) = \frac{1}{2\pi} \int_{\R^2_+} \big ( \frac{ (x-y)^\bot}{|x-y|^2} -\frac{(x-y^*)^\bot}{|x-y^*|^2} \big ) \omega (t,y) \dd y.\label{eq.lem.biot-savart}
\end{equation}
Here $x^\bot = (-x_2,x_1)^\top $ and $y^* = (y_1,-y_2)^\top $.
\end{lem}

\noindent {\it Proof of Theorem \ref{thm.intro.1}.} We give a proof of Theorem \ref{thm.intro.1} by admitting Lemmas \ref{lem.spatial.decay.vorticity} and \ref{lem.biot-savart}. The proofs of  these lemmas will be postponed to the latter sections. From \eqref{est.lem.spatial.decay.vorticity} and \eqref{eq.lem.biot-savart} we observe that
\begin{align}
0=\lim_{x_2 \downarrow 0} u_1 (t,x) = \frac{1}{\pi} \int_{\R^2+} \frac{y_2}{(x_1-y_1)^2 + y_2^2} \omega (t,y) \dd y,\label{proof.thm.main.1}
\end{align}
by the Lebesgue convergence theorem. Then {\rm (C4)} implies that the integrand of the right-hand side of \eqref{proof.thm.main.1} has to be zero, that is, $\omega (t,x) =0$ in $(-\infty,0)\times \R^2_+$. Then for all $t$ the velocity $u(t)$ is harmonic and bounded in $\R^2_+$ and vanishes on $\partial\R^2_+$. Hence $u$ must be zero by the classical Liouville theorem for harmonic functions. The proof is now complete.

\subsection{Temporal decay of vorticity}\label{sec.temporal.decay}

The main result of this section is the following lemma.

\begin{lem}\label{lem.temporal.decay.vorticity} Under the conditions {\rm (C1)}, {\rm (C2)}, and {\rm (C3)} of Theorem \ref{thm.intro.1} the vorticity $\omega$ satisfies
\begin{align}
\| \omega (t) \|_\infty & \leq C (-t)^{-1} |\log (-t)|^2,~~~~~~~~-\infty < t < -2,\label{est.lem.temporal.decay.vorticity.1}\\
\| (-\partial_1^2 )^\frac12 \omega (t) \|_\infty & \leq C (-t)^{-\frac32} |\log (-t )|^4,~~~~~~~~-\infty < t < -2.\label{est.lem.temporal.decay.vorticity.2}
\end{align}

\end{lem}

Lemma \ref{lem.temporal.decay.vorticity} is proved by estimating the integral equations for the vorticity field in Proposition \ref{prop.vorticity.equation}. To this end we first establish the $L^\infty - L^\infty$ estimates for the operators in \eqref{eq.lem.integral}.

\begin{lem}\label{lem.semigroup.estimate} Assume that $v \in (C (\overline{\R^2_+}))^2$ with $v =0$ on $\partial\R^2_+$ and $g\in BC(\R)$. Then
\begin{align}
\| (-\partial_1^2 )^{\frac{1}{2}} T(t) v  \|_\infty & \leq C t^{-1} \| v \|_\infty,\label{est.lem.semigroup.estimate.0}\\
\| (-\partial_1^2 )^{\frac{k}{2}} T(t) \partial_i v  \|_\infty & \leq C t^{-1-\frac{k}{2}} \| v \|_\infty,~~~~~~~k=0,1,~i=1,2,\label{est.lem.semigroup.estimate.1}\\
\| (-\partial_1^2 )^{\frac{k}{2}} e^{t B} ( \partial_1^l g \delta_{\partial\R^2_+} ) \|_\infty & \leq C t^{-\frac{1+k+l}{2}} \| g\|_\infty,~~~~~k,l=0,1.\label{est.lem.semigroup.estimate.2}
\end{align}
Moreover, if $F=u\otimes u$ with $u\in (C^2 (\overline{\R^2_+}))^2$ satisfying ${\rm div}\, u=0$ in $\R^2_+$ and $u=0$ on $\partial\R^2_+$ then
\begin{align}
\| (-\partial_1^2)^\frac12 T(t) {\rm div} F \|_\infty \leq C \| u \|_\infty \min\{ t^{-1}  \| \omega \|_\infty,~ t^{-\frac12} \|\nabla \omega \|_\infty\} .\label{est.lem.semigroup.estimate.3}
\end{align}
\end{lem}

\noindent {\it Proof.} As in \cite[Proposition 5.1]{Mae}, using the Fourier transform,  we can derive the pointwise estimate for the kernel $K(t,x,y)$ of $e^{tB}$ such as
\begin{align}
& ~~~ |(-\partial_1^2)^\frac{k}{2} \partial_1^l \partial_2^j  K (t,x,y) | \nonumber \\
& \leq  C t^{-\frac{k+l+2}{2}} \bigg (1+\frac{|(x_1-y_1)/\sqrt{t}|^{2+k+l}}{\log (e+|(x_1-y_1)/\sqrt{t}|^2) }  + |(x_2-y_2)/\sqrt{t}|^{2+k+l+j} \bigg )^{-1}.\label{est.prop.pointwise.kernel}
\end{align}
Since $e^{t B}$ commutes with $\partial_1$, the estimates \eqref{est.lem.semigroup.estimate.0} and \eqref{est.lem.semigroup.estimate.2} are immediate from \eqref{est.prop.pointwise.kernel}. As for \eqref{est.lem.semigroup.estimate.1}, we give a proof only for the case $k=1$ and $i=2$. The other cases are proved in the same manner.  By the definition of $T(t)$ in \eqref{def.T(t)} we have
\begin{align*}
\langle (-\partial_1^2)^\frac{1}{2}T(t) \partial_2 v, f\rangle _{L^2} & = \langle \partial_2 v_1, \partial_2 (-\partial_1^2)^\frac{1}{2} e^{t B} f \rangle _{L^2}   - \langle \partial_2 v_2, \partial_1 (-\partial_1^2)^\frac{1}{2} e^{t B} f\rangle _{L^2}\\
& = - \langle v_1, \partial_2^2 (-\partial_1^2)^\frac{1}{2} e^{t B} f \rangle _{L^2}   + \langle v_2, \partial_1 \partial_2 (-\partial_1^2)^\frac{1}{2} e^{t B} f\rangle _{L^2}.
\end{align*}
Here we have used the integration by parts and the boundary condition $v=0$ on $\partial\R^2_+$. Since \eqref{est.prop.pointwise.kernel} implies $\|\partial_2^2 (-\partial_1^2)^\frac{1}{2} e^{t B} f\|_1 + \| \partial_2^2 (-\partial_1^2)^\frac{1}{2} e^{t B} f\|_{L^1}\leq C t^{-3/2} \| f\|_{L^1}$, we obtain \eqref{est.lem.semigroup.estimate.2} by the duality argument. Finally we show \eqref{est.lem.semigroup.estimate.3}. Set $v={\rm div}~F$. Note that $v$ vanishes on the boundary by the assumption. Then again by the definition of $T(t)$ we have
\begin{align*}
\langle (-\partial_1^2)^\frac{1}{2}T(t) {\rm div}\, F , f\rangle _{L^2} & = \langle v_1, \partial_2 (-\partial_1^2)^\frac{1}{2} e^{t B} f \rangle _{L^2}   - \langle v_2, \partial_1 (-\partial_1^2)^\frac{1}{2} e^{t B} f\rangle _{L^2}\\
& = \langle - \nabla^\bot \cdot v, (-\partial_1^2)^\frac{1}{2} e^{t B} f \rangle _{L^2}\\
& = \langle u\cdot \nabla \omega ,  (-\partial_1^2)^\frac{1}{2} e^{t B} f \rangle _{L^2} = - \langle u \omega , \nabla  (-\partial_1^2)^\frac{1}{2} e^{t B} f \rangle _{L^2}.
\end{align*}
Here we have used the equality $-\nabla^\bot \cdot {\rm div}\, (u\otimes u) = u\cdot \nabla \omega = \nabla \cdot (u \omega )$. By using the estimates  $\| \nabla ^k (-\partial_1^2)^\frac{1}{2} e^{t B} f\|_{L^1} \leq C t^{-(1+k)/2} \| f\|_{L^1}$ for $k=0,1$, we obtain \eqref{est.lem.semigroup.estimate.3} by the duality argument. The proof is now complete.

\vspace{0.5cm}

\noindent {\it Proof of Lemma \ref{lem.temporal.decay.vorticity}.} By Proposition \ref{prop.vorticity.equation} the vorticity $\omega$ satisfies the integral equation \eqref{eq.lem.integral} for $-\infty<s<2 t< t<0$ with $F (t) =-u (t) \otimes u (t)$. We set
\begin{align*}
I(t,s) &= T(t-s) u(s),~~~~~~~~~~~~~~~ II(t,s) = \int_s^t T(t-\tau ){\rm div} \, F(\tau ) \dd \tau,\\
III(t,s) &= \int_s^t e^{(t-\tau )B} (\partial_1 p_{F (\tau)} \delta_{\partial\R^2_+} )  \dd\tau.
\end{align*}
For $I(t,s)$ we have from \eqref{est.T.1}, \eqref{est.lem.semigroup.estimate.0},  and (C3) that
\begin{align}
\|(-\partial_1^2)^\frac{k}{2} I (t,s) \|_\infty\leq C (t-s)^{-\frac{1+k}{2}} \| u(s) \|_\infty\leq C (t-s)^{-\frac{1+k}{2}}  (-s)^{-\frac{1}{2}}\rightarrow 0~~~{\rm  as}~~ s\rightarrow -\infty.\label{proof.lem.temporal.decay.vorticity.1}
\end{align}
Next we consider the term $II(t,s)$. When $\tau <t-1/t^2$ we apply \eqref{est.lem.semigroup.estimate.1} and get
\[
\|T(t-\tau) {\rm div} \, F(\tau) \|_\infty\leq C (t-\tau)^{-1} \| u (\tau )\|_\infty^2 \leq C(t-\tau)^{-1} (-\tau )^{-1}
\]
by (C3), and we also have from \eqref{est.lem.semigroup.estimate.3} and (C3) that  
\[
\| (-\partial_1^2)^{\frac12} T(t-\tau) {\rm div}\, F(\tau) \|_\infty\leq C (t-\tau)^{-1} (-\tau )^{-\frac12} \| \omega (\tau ) \|_\infty.
\]
When $t-1/t^2\leq \tau <t$ we use \eqref{est.T.1}, \eqref{est.lem.semigroup.estimate.0}, and \eqref{est.lem.semigroup.estimate.3} to get
\[
\| T(t-\tau) {\rm div}\, F(\tau) \|_\infty\leq C (t-\tau)^{-\frac12} \| u (\tau ) \cdot \nabla u (\tau )\|_\infty\leq C (t-\tau)^{-\frac12} (-\tau )^{-\frac12}\]
and
\[
\| (-\partial_1^2)^{\frac12} T(t-\tau) {\rm div}\, F(\tau) \|_\infty\leq C (t-\tau)^{-\frac12} \|u (\tau ) \|_\infty \| \nabla \omega  (\tau )\|_\infty\leq C (t-\tau)^{-\frac12} (-\tau )^{-\frac12}.
\]
Collecting these, for $t<-2$ we have arrived at
\begin{align}
\lim_{s\rightarrow -\infty} \| II(t,s) \|_\infty & \leq C \int_{-\infty}^{t-\frac{1}{t^2}} (t-\tau)^{-1} (-\tau )^{-1} \dd\tau + C \int_{t-\frac{1}{t^2}}^t (t-\tau)^{-\frac12} (-\tau )^{-\frac12} \dd\tau \nonumber \\
& \leq C (-t)^{-1} \log (-t) ,\label{proof.lem.temporal.decay.vorticity.2}
\end{align}
and
\begin{align}
\lim_{s\rightarrow -\infty}\| (-\partial_1^2)^\frac12  II(t,s) \|_\infty %&\leq C \int_{-\infty}^{t-\frac{1}{t^2}} (t-\tau)^{-1} (-\tau )^{-\frac12} \|\omega (\tau ) \|_\infty  \dd\tau \nonumber \\
%& ~~~ + C \int_{t-\frac{1}{t^2}}^t (t-\tau)^{-\frac12} (-\tau )^{-\frac12} \dd\tau \nonumber \\
& \leq C \int_{-\infty}^{t-\frac{1}{t^2}} (t-\tau)^{-1} (-\tau )^{-\frac12} \|\omega (\tau ) \|_\infty  \dd\tau + C (-t)^{-\frac32}.\label{proof.lem.temporal.decay.vorticity.3}
\end{align}
Finally we estimate $III(t,s)$. To this end we recall that $p_{F}$ is the restriction of the function $-{\rm div}\, {\rm div} \, (-\Delta_{\R^2})^{-1} \tilde F$ on $\R^2_+$; see the proof of Proposition  \ref{prop.poisson} for details and the definition of $\tilde F$. Then we decompose $\partial_1 p_{F(\tau)}$ as
\begin{align*}
\partial_1 p_{F(\tau )} = \sum_{j=1}^3 \partial_1 p_{F(\tau),j} = -\big (  \int_0^{\frac{1}{\tau^4}} + \int_{\frac{1}{\tau^4}}^{\tau^4} + \int_{\tau^4} ^\infty \big )  \partial_1 {\rm div}\, {\rm div}\,  G (\theta ) * \tilde F (\tau) \dd\theta.
\end{align*}
Here $G(\theta,x)$ is the two-dimensional Gaussian. Firstly we observe that
\begin{align*}
\| (-\partial_1^2)^{\frac{k}{2}} e^{(t-\tau)B} (\partial_1 p_{F(\tau),1} \delta_{\partial\R^2_+}) \|_\infty &  \leq C (t-\tau)^{-\frac12 - \frac{k \kappa}{2}} \int_0^{\frac{1}{\tau^4}}\theta^{-\frac12 -\frac{k(1-\kappa)}{2}}  \| {\rm div}\, {\rm div}\, \tilde F(\tau )\|_\infty \dd\tau \\
& \leq C  (t-\tau)^{-\frac12 - \frac{k\kappa}{2}} (-\tau)^{-2 +2 k (1- \kappa )}
\end{align*}
for $k=0,1$ and $\kappa\in (0,1)$, where we have applied  \eqref{est.lem.semigroup.estimate.2} and the interpolation argument using $(-\partial_1^2)^{1/2} = (-\partial_1^2)^{\kappa/2} (-\partial_1^2)^{(1-\kappa)/2}$ when $k=1$. By taking $\kappa$ close to $1$ we thus obtain
\begin{align}
\| \int_{-\infty}^t (-\partial_1^2)^{\frac{k}{2}} e^{(t-\tau)B} (\partial_1 p_{F(\tau),1} \delta_{\partial\R^2_+}) \dd\tau \|_\infty \leq C (-t)^{-\frac32},~~~~~~~~~~k=0,1,~-\infty<t<-2.\label{proof.lem.temporal.decay.vorticity.4}
\end{align}
The estimate of $\partial_1 p_{F(\tau),3}$ is easily calculated as
\[
\| (-\partial_1^2)^{\frac{k}{2}} \partial_1 p_{F(\tau),3} \|_\infty \leq C  \int_{\tau^4}^\infty \theta^{-\frac{3+k}{2}} \dd\theta \| \tilde F(\tau) \|_\infty\leq C (-\tau)^{-3},~~~~~k=0,1.
\]
Hence we have from  \eqref{est.lem.semigroup.estimate.2},
\begin{align}
\| \int_{-\infty}^t (-\partial_1^2)^{\frac{k}{2}} e^{(t-\tau)B} (\partial_1 p_{F(\tau),3} \delta_{\partial\R^2_+}) \dd\tau \|_\infty \leq C (-t)^{-\frac32},~~~~~~~~~~k=0,1,~-\infty<t<-2.\label{proof.lem.temporal.decay.vorticity.5}
\end{align}
Now we consider the term related with $\partial_1 p_{F(\tau),2}$. By the definition of $\tilde F$ in Proposition \ref{prop.poisson} we take the even extension for $u_1$ and the odd extension for $u_2$. Each extension is denoted by $\tilde u_i$. This extension leads to the odd extension $\tilde \omega$ of the vorticity $\omega$. Then it is straightforward to see ${\rm div} \, \tilde F = - \tilde u^\bot \tilde  \omega - \nabla | \tilde u|^2/2$ with $\tilde u^\bot = (-\tilde u_2,\tilde u_1)^\top$, and thus, ${\rm div}\, {\rm div}\, \tilde F = -{\rm div}\, (\tilde u^\bot \tilde \omega ) - \Delta |\tilde u |^2/2$.  Hence we have
\begin{align*}
\partial_1 p_{F(\tau),2} & = \int_{\frac{1}{\tau^4}}^{\tau^4} \partial_1 {\div} \, G(\theta )  \dd\theta * (\tilde u^\bot \tilde \omega ) (\tau ) + \frac12 \int_{\frac{1}{\tau^4}}^{\tau^4}\partial_1 \Delta G(\theta ) \dd\theta * |\tilde u |^2   (\tau )  \nonumber  \\
& =  \int_{\frac{1}{\tau^4}}^{\tau^4} \partial_1 {\div}\, G(\theta ) \dd\theta  * (\tilde u^\bot \tilde \omega )  (\tau )  + \frac12 \partial_1 G(\tau^4) * |\tilde u|^2 - \frac12  G(\tau^{-4}) * \partial_1|\tilde u|^2 (\tau ) . 
\end{align*}
Since $\partial_1 |\tilde u|^2 =0$ on $\partial\R^2_+$ we have $\| G(\tau^{-4}) * \partial_1|\tilde u|^2 \|_{L^\infty (\partial\R^2_+)}\leq C (-\tau)^{-2}\| \partial_1 |\tilde u |^2 \|_{C^1} \leq C (-\tau)^{-2}$. Hence it follows that
\begin{align}
\| \partial_1 p_{F(\tau),2}\|_{L^\infty (\partial\R^2_+)} & \leq C \| \tilde u(\tau ) \|_\infty \| \tilde \omega (\tau ) \|_\infty \log (-\tau ) + C (-\tau )^{-2} \nonumber \\
& \leq C  \| \omega (\tau ) \|_\infty (-\tau)^{-\frac12} \log (-\tau )+ C (-\tau )^{-2}.\label{proof.lem.temporal.decay.vorticity.6}
\end{align}
When $\tau<t-1/t^4$ we have from \eqref{est.lem.semigroup.estimate.2} that
\begin{align*}
\| e^{(t-\tau)B} (\partial_1 p_{F(\tau),2} \delta_{\partial\R^2_+}) \|_\infty & \leq C (t-\tau)^{-1} \int_{\frac{1}{\tau^4}}^{\tau^4} \theta^{-1} \| \tilde F(\tau)\|_\infty \dd\theta \leq C  (t-\tau)^{-1} (-\tau )^{-1} \log (-\tau ),
\end{align*}
while  \eqref{proof.lem.temporal.decay.vorticity.6} implies
\begin{align*}
\| (-\partial_1^2)^\frac12 e^{(t-\tau)B} (\partial_1 p_{F(\tau),2} \delta_{\partial\R^2_+}) \|_\infty & \leq C (t-\tau)^{-1} \big ( \| \omega (\tau ) \|_\infty (-\tau)^{-\frac12}  \log (-\tau )+ (-\tau )^{-2} \big ).
\end{align*}
As for the case $t-1/t^4 \leq \tau <t$, we have for $k=0,1$,
\begin{align*}
\| (-\partial_1^2)^{\frac{k}{2}} e^{(t-\tau)B} (\partial_1 p_{F(\tau),2} \delta_{\partial\R^2_+}) \|_\infty & \leq C (t-\tau)^{-\frac12} \int_{\frac{1}{\tau^4}}^{\tau^4} \theta^{-1} \| \tilde F(\tau)\|_{C^2}\dd\theta \\
& \leq C  (t-\tau)^{-\frac12} \log (-\tau ).
\end{align*}
Combining the above three yields
\begin{align}
\| \int_{-\infty}^t  e^{(t-\tau)B} (\partial_1 p_{F(\tau),2} \delta_{\partial\R^2_+}) \dd\tau \|_\infty \leq C (-t)^{-1} |\log (-t)|^2,~~~~~~-\infty <t<-2,\label{proof.lem.temporal.decay.vorticity.7}
\end{align}
and
\begin{align}
& ~~~ \| (-\partial_1^2)^\frac12 \int_{-\infty}^t  e^{(t-\tau)B} (\partial_1 p_{F(\tau),2} \delta_{\partial\R^2_+}) \dd\tau \|_\infty \nonumber \\
& \leq C \int_{-\infty}^{t-1/t^4} (t-\tau )^{-1} \|\omega (\tau)\|_\infty (-\tau )^{-\frac12} \log (-\tau ) \dd\tau + C (-t)^{-\frac32} ,~~~~~~-\infty <t<-2.\label{proof.lem.temporal.decay.vorticity.8}
\end{align}
The estimates \eqref{proof.lem.temporal.decay.vorticity.1}, \eqref{proof.lem.temporal.decay.vorticity.2}, \eqref{proof.lem.temporal.decay.vorticity.4}, \eqref{proof.lem.temporal.decay.vorticity.5}, \eqref{proof.lem.temporal.decay.vorticity.7} imply \eqref{est.lem.temporal.decay.vorticity.1}, and the estimates  \eqref{proof.lem.temporal.decay.vorticity.1}, \eqref{proof.lem.temporal.decay.vorticity.3}, \eqref{proof.lem.temporal.decay.vorticity.4}, \eqref{proof.lem.temporal.decay.vorticity.5}, \eqref{proof.lem.temporal.decay.vorticity.8} together with
\eqref{est.lem.temporal.decay.vorticity.1} give \eqref{est.lem.temporal.decay.vorticity.2}. The proof is complete.

\begin{rem}\label{rem.lem.temporal.decay.vorticity}{\rm The proof of Lemma \ref{lem.temporal.decay.vorticity} implies that, from \eqref{est.lem.temporal.decay.vorticity.1} and \eqref{proof.lem.temporal.decay.vorticity.6},
\[
\| \partial_1 p_{F(\tau),2} \|_{L^\infty (\partial\R^2_+)} \leq C (-\tau)^{-\frac32} |\log (-\tau)|^3,~~~~~-\infty <\tau<-2.
\]
Since it is easy to see  $\| \partial_1 p_{F(\tau),i} \|_{L^\infty (\partial\R^2_+)} \leq (-\tau )^{-3/2}$ for $i=1,3$, we have
\begin{align}
 \| \partial_1 p_{F(\tau)} \|_{L^\infty (\partial\R^2_+)} \leq C (-\tau)^{-\frac32} |\log (-\tau)|^3,~~~~~~-\infty < \tau <-2.\label{est.rem.lem.temporal.decay.vorticity.1}
\end{align}
This estimate will be used later.

}
\end{rem}

\subsection{Spatial decay of vorticity - proof of Lemma \ref{lem.spatial.decay.vorticity}}\label{sec.spatial.decay}

In this section we derive spatial decay  of the vorticity field and complete the proof of Lemma \ref{lem.spatial.decay.vorticity}. The key idea is to regard \eqref{eq.vorticity.boundary} as the Neumann boundary condition $\partial_2\omega =g$  with the inhomogeneous term $g=-(-\partial_1^2)^{1/2} \omega |_{x_2=0} - \partial_1 p_F$. Then we use a representation formula of the vorticity in terms of the fundamental solution for the heat-transport operator $\partial_t - \Delta + \tilde u \cdot \nabla$ in $\R^2$, whose precise pointwise estimate has already been established by \cite{CarLo, MaTo}. Here $\tilde u=(\tilde u_1, \tilde u_2)^\top$ is the extension of $u$ to $\R^2$, where $\tilde u_1$, $\tilde u_2$ are the even, odd extensions of $u_1,~u_2$, respectively. Note that this extension preserves the divergence-free condition when $u_2=0$ on $\partial\R^2_+$. The scaling invariant assumption (C3) is essential in establishing the spatial decay of the vorticity, for it leads to the global Gaussian estimate for the fundamental solution.  We start from the following lemma.

\begin{lem}\label{lem.heat-transport} Under the conditions {\rm (C1)}, {\rm (C2)}, and {\rm (C3)} of Theorem \ref{thm.intro.1} the vorticity $\omega$ is expressed as
\begin{align}
\omega (t) = \Gamma_u (t,s) \omega (s) -  \int_s^t \Gamma_u (t,\tau ) ( g (\tau )  \delta_{\partial\R^2_+} )  \dd\tau,~~~~~~-\infty <s<t<0,\label{eq.lem.heat-transport}
\end{align}
with $g (\tau ) =-(-\partial_1^2)^{1/2} \omega  (\tau ) |_{x_2=0} - \partial_1 p_{F(\tau )}$. Here $\Gamma_u(t,s)$ is the evolution operator  defined by
\[
\Gamma_u (t,s)f =\int_{\R^2_+} \big ( \Gamma_{\tilde u} (t,x;s,y) + \Gamma_{\tilde u} (t,x;s,y^*) \big ) f (y) \dd y,
\]
where $\Gamma_{\tilde u} (t,x;s,y)$ is the fundamental solution to  the heat-transport equations
\begin{align}
\partial_t w -\Delta w + \tilde u \cdot \nabla w =0~~~~~~~~~{\rm in} ~~ (-\infty,0) \times \R^2. \label{eq.heat-transport.whole}
\end{align}
Moreover, it follows that
\begin{align}
\| \Gamma_u (t,s) f\|_p & \leq C (t-s)^{-\frac{1}{q}+\frac{1}{p}} \| f \|_q,~~~~~~-\infty<s<t<0,~~~1\leq q\leq p\leq \infty, \label{est.lem.heat-transport.1}\\
0<\Gamma_{\tilde u} (t,x;s,y) & \leq C_1 (t-s)^{-1} \exp \big ( -C_2\frac{|x-y|^2}{t-s} \big ).\label{est.lem.heat-transport.2}
\end{align}
Here $C_1$ and $C_2$ depend only on $M= \displaystyle \sup_{-\infty < t <0} (-t)^{1/2} \| u(t) \|_\infty$.

\end{lem}

\begin{rem}{\rm In \eqref{eq.lem.heat-transport} the term $\Gamma_u (t,\tau ) ( g (\tau )  \delta_{\partial\R^2_+} )$ is defined as 
\begin{align}
\Gamma_u (t,\tau ) ( g (\tau )  \delta_{\partial\R^2_+} ) (x) =2  \int_\R \Gamma_{\tilde u} (t,x;\tau ,y_1,0) g (\tau , y_1 ) \dd y_1.\label{rem.lem.heat-transport}
\end{align}

}
\end{rem}

\noindent {\it Proof of Lemma \ref{lem.heat-transport}.} The existence of fundamental solutions to \eqref{eq.heat-transport.whole} is classical under the assumption of (C1); cf. \cite{GiGiSa}. The estimate \eqref{est.lem.heat-transport.1} is a consequence of \cite[Theorem 1]{CarLo} and the definition of $\Gamma_u(t,s)$. As for \eqref{est.lem.heat-transport.2}, we have from \cite[Theorem 3]{CarLo} that
\begin{align}
\Gamma_{\tilde u} (t,x;s,y) \leq \frac{1}{4\pi (t-s)} \exp \bigg ( - \frac{1}{4 (t-s)}  \big ( |x-y| - \int_s^t \| u(\tau ) \|_\infty \dd \tau \big )_+^2 \bigg ). \label{proof.lem.heat-transport.1}
\end{align}
Here $(\alpha)_+=\max\{0,\alpha\}$ for $\alpha\in \R$. The condition (C3) yields
\[
\int_s^t \| u(\tau ) \|_\infty \dd\tau \leq M \int_s^t (-\tau)^{-\frac12} \dd\tau \leq 2 M|t-s|^{\frac12},~~~~~ M= \sup_{-\infty < t <0} (-t)^{\frac12} \| u(t) \|_\infty.
\]
Hence if $|x-y|\geq 4 (t-s)^{1/2}$ then  \eqref{proof.lem.heat-transport.1} implies \eqref{est.lem.heat-transport.2}. On the other hand, if $|x-y|\leq 4 M (t-s)^{1/2}$ then again from \eqref{proof.lem.heat-transport.1}  we have
\begin{align*}
\Gamma_{\tilde u} (t,x;s,y) \leq \frac{1}{4\pi (t-s)} = \frac{1}{4\pi (t-s)} e^{\frac{|x-y|^2}{t-s}} e^{-\frac{|x-y|^2}{t-s}} \leq \frac{e^{16 M^2}}{4\pi (t-s)} e^{-\frac{|x-y|^2}{t-s}},
\end{align*}
which is the desired estimate. The positivity of $\Gamma_{\tilde u} (t,x;s,y)$ is a consequence of the strong maximal principle and the details are omitted here. The representation \eqref{eq.lem.heat-transport} is derived from the fact that the equation $\partial_t\omega - \Delta \omega + u \cdot \nabla \omega =0$ in $(-\infty,0)\times \R^2_+$ with the Neumann boundary condition $\partial_2 \omega = g$ on $\partial\R^2_+$ is equivalent with the equation
\begin{equation}
\partial_t \tilde w -\Delta \tilde w + \tilde u\cdot \nabla \tilde w =-2 g \delta_{\partial\R^2_+}~~~{\rm in}~(-\infty,0)\times \R^2,
\end{equation}
where $\tilde w$ is the even extension of $\omega$ to $\R^2$. The proof is complete.

\vspace{0.5cm}

\noindent {\it Proof of Lemma \ref{lem.spatial.decay.vorticity}.} By (C1), \eqref{est.lem.temporal.decay.vorticity.2}, and  \eqref{est.rem.lem.temporal.decay.vorticity.1}  the function $g(t) = -(-\partial_1^2)^{1/2} \omega  (t) |_{x_2=0} - \partial_1 p_{F(t)}$ is estimated as
\begin{align}
\| g(t) \|_{L^\infty (\partial\R^2_+)} \leq C (-t)^{-\frac{3-\epsilon}{2}}~~~~~~~-\infty < t < 0,~\epsilon\in (0,1). \label{proof.lem.spatial.decay.vorticity.1}
\end{align}
The estimate \eqref{est.lem.heat-transport.2} and the representation \eqref{rem.lem.heat-transport} lead to
\[
\| \Gamma_u (t,\tau ) (g (\tau ) \delta_{\partial\R^2_+} )\|_\infty \leq C (t-\tau )^{-\frac12} (-\tau )^{-\frac{3-\epsilon}{2}}
\]
for $\tau <0$ and $\epsilon\in (0,1)$. On the other hand, we have from \eqref{est.lem.temporal.decay.vorticity.1} and  \eqref{est.lem.heat-transport.1} that $\| \Gamma_u (t,s) \omega (s) \|_\infty \leq C (-s)^{-1}|\log (-s)|^2$ for $s\ll -1$. Thus by taking the limit $s\rightarrow -\infty$ in \eqref{eq.lem.heat-transport} we arrive at the expression
\begin{align}
\omega (t,x) = - 2 \int_{-\infty}^t \int_\R \Gamma_u (t,x;\tau, y_1,0) g (\tau, y_1 ) \dd y_1 \dd\tau,~~~~~~~t<0,~x\in \R^2_+.\label{proof.lem.spatial.decay.vorticity.2}
\end{align}
Let $\theta\in (0,1-\epsilon)$. Then from \eqref{est.lem.heat-transport.2} and \eqref{proof.lem.spatial.decay.vorticity.1} we have
\begin{align}
x_2^{1+\theta} |\omega (t,x)| & \leq C\int_{-\infty}^t \int_\R (t-\tau)^{-1+\frac{1+\theta}{2}} e^{-c\frac{(x_1-y_1)^2}{t-\tau}} | g (\tau, y_1 ) | \dd y_1 \nonumber \\
& \leq C\int_{-\infty}^t  (t-\tau )^{\frac{\theta}{2}} (-\tau )^{-\frac{3-\epsilon}{2}} \dd\tau  \leq C (-t )^{-\frac{1-\theta-\epsilon}{2}}.
\end{align}
It is easy to see that the same argument with (C1) also yields $\displaystyle \sup_{-1<t<0,x\in \R^2_+} x_2^{1+\theta} |\omega (t,x)|<\infty$. The proof is complete.

\subsection{Representation of solutions by the Biot-Savart law}\label{sec.bs.law}

In this section we give a proof of Lemma \ref{lem.biot-savart}. To this end we denote by $v(t,x)$ the right-hand side of \eqref{eq.lem.biot-savart}, which is well-defined by \eqref{est.lem.spatial.decay.vorticity} and the estimate
\[
\int_{\R^2_+} \big | \frac{(x-y)^\bot}{|x-y|^2} - \frac{(x-y^*)^\bot}{|x-y^*|^{2}} \big |  ~ (1+ y_2)^{-1-\theta} \dd y \leq C\int_{\R^2_+} \frac{y_2}{|x-y| |x-y^*|} (1+y_2)^{-1-\theta} \dd y <\infty.
\]
In particular, $v$ is uniformly  bounded in $(-\infty,0)\times \R^2_+$. The goal is thus to show $u=v$.  Since both $u$ and $v$ satisfy the divergence-free condition and their vorticity fields are given by the same $\omega$, the difference $w=u-v$ is harmonic in $\R^2_+$. Moreover, $u$ and $v_2$ vanishes on the boundary by the no-slip boundary condition and the definition of $v$. Hence, due to the Liouville theorem for harmonic functions we only need to prove the fact $v_1=0$ on $\partial\R^2_+$. We first note that $v_1$ is written as
\begin{align}
v_1 (t,\cdot, x_2) & = \partial_2 \int_0^{x_2} e^{-(x_2-y_2)(-\partial_1^2)^\frac12}\int_{y_2}^\infty e^{-(z_2-y_2)(-\partial_1^2)^\frac12} \omega (t,\cdot, z_2) \dd z_2 \dd y_2 \nonumber \\
& = \int_{x_2}^\infty e^{-(y_2-x_2)(-\partial_1^2)^\frac12}\omega (t,\cdot,y_2) \dd y_2  \nonumber \\
& ~~~ -  \int_0^{x_2}\int_{y_2}^\infty (-\partial_1^2)^{\frac12} e^{-(x_2 - 2 y_2+z_2) (-\partial_1^2)^\frac12} \omega (t,\cdot,z_2) \dd z_2 \dd y_2. \label{proof.lem.biot-savart.1}
\end{align}
The last term of the right-hand side of \eqref{proof.lem.biot-savart.1} vanishes on $\partial\R^2_+$, so we focus on the first term which we will denote by $v_{1,1} (t,x)$. Fix any $\delta>0$ and let $-t> 2 \delta$ and $x_2>\delta$. For sufficiently small $\epsilon\in (0, \delta/4)$ we denote by $\omega_\epsilon (t,x) =\int_{-\infty}^\delta  \int_{\R^2_+} \eta_\epsilon (t-s, x-y) \omega (s,y ) \dd y \dd s$ the mollification of $\omega$. The mollifier $\eta_\epsilon$ is taken so that ${\rm supp} ~ \eta_\epsilon \subset \{ (t,x)\in \R^{3} ~|~|t|^2 + |x|^2 <\epsilon^2\}$. Then $\omega_\epsilon$ satisfies
\begin{align}
\partial_t \omega_{\epsilon} (t,x) & =  \Delta \omega_\epsilon (t,x) - \nabla \cdot (u \omega )_\epsilon (t,x) + F_\epsilon (t,x),\label{proof.lem.biot-savart.2}
\end{align}
where
\begin{align*}
(u\omega)_\epsilon (t,x) & =  \int_{-\infty}^\delta \int_{\R^2_+} \eta_\epsilon (t-s,x-y) u \omega (s,y) \dd y\dd s,\\
F_\epsilon (t,x) & =  -\eta_\epsilon (t-\delta) *\omega (\delta)(x) \\
& ~~ -\int_{-\infty}^\delta \int_{\partial \R^2_+} \big ( \eta_\epsilon  (t-s, x-y) \partial_2 \omega (s,y) + \partial_2 \eta_\epsilon (t-s,x-y) \omega (s,y ) \big ) \dd y_1 \dd s.
 \end{align*}
By \eqref{est.lem.spatial.decay.vorticity} and the definition of $\eta_\epsilon$ each term in \eqref{proof.lem.biot-savart.2} has the same spatial decay as $\omega$. Set $v_{1,1,\epsilon} (t,\cdot,x_2) = \int_{x_2}^\infty e^{-(y_2-x_2)(-\partial_1^2)^\frac12}\omega_\epsilon (t,\cdot,y_2) \dd y_2$.
Then we verify the calculation
\begin{align*}
\partial_t v_{1,1,\epsilon} (t,\cdot,x_2) = \int_{x_2}^\infty e^{-(y_2-x_2)(-\partial_1^2)^\frac12} \big ( \Delta \omega_\epsilon - \nabla \cdot (u\omega)_\epsilon + F_\epsilon \big ) (t,\cdot,y_2)  \dd y_2,
\end{align*}
and the integration by pars yields
\begin{align}
\partial_t v_{1,1,\epsilon} (t,\cdot,x_2) & = -\partial_2 \omega_\epsilon (t,\cdot,x_2) - (-\partial_1^2)^\frac12 \omega_\epsilon (t,\cdot,x_2) \nonumber \\
& ~~~ -   \int_{x_2}^\infty \nabla_x \cdot e^{-(y_2-x_2)(-\partial_1^2)^\frac12} (u\omega )_\epsilon (t,\cdot,y_2) \dd y_2 \nonumber \\
& ~~~~~~ + (u_2 \omega )_\epsilon (t,\cdot,x_2) +  \int_{x_2}^\infty e^{-(y_2-x_2)(-\partial_1^2)^\frac12} F_\epsilon (t,\cdot y_2) \dd y_2.\label{proof.lem.biot-savart.3}
\end{align}
From (C1) and \eqref{est.lem.spatial.decay.vorticity} it is easy to see that the following convergence holds in the limit $\epsilon\rightarrow 0$ uniformly on each compact set of $\{(t,x)~|~t<-2\delta, x_2>\delta\}$:
\begin{align*}
-\partial_2 \omega_\epsilon (t) - (-\partial_1^2)^\frac12 \omega_\epsilon (t)+ (u_2 \omega )_\epsilon (t) & ~\rightarrow ~ -\partial_2 \omega (t) - (-\partial_1^2)^\frac12 \omega (t)+ u_2 \omega  (t),\\
\int_{x_2}^\infty \nabla\cdot e^{-(y_2-x_2)(-\partial_1^2)^\frac12} (u\omega )_\epsilon (t,\cdot,y_2) \dd y_2  & ~\rightarrow ~ \int_{x_2}^\infty \nabla_x \cdot e^{-(y_2-x_2)(-\partial_1^2)^\frac12} (u\omega ) (t,\cdot,y_2) \dd y_2, \\
\int_{x_2}^\infty e^{-(y_2-x_2)(-\partial_1^2)^\frac12} F_\epsilon (t,\cdot,y_2) \dd y_2 & ~ \rightarrow ~ 0.
\end{align*}
Thus we have for $s<t<-2\delta$ and $x_2>\delta$,
\begin{align}
v_{1,1} (t) - v_{1,1} (s) & = \int_{s}^{t} \big ( -\partial_2 \omega (\tau ) - (-\partial_1^2)^\frac12 \omega (\tau )+ u_2 \omega  (\tau ) \big ) \dd\tau\nonumber \\
& ~~~ - \int_{s}^{t}  \int_{x_2}^\infty \nabla_x \cdot e^{-(y_2-x_2)(-\partial_1^2)^\frac12} (u\omega ) (\tau,\cdot,y_2 ) \dd y_2 \dd \tau.\label{proof.lem.biot-savart.4}
\end{align}
Since $\delta>0$ is arbitrary we may take $x_2\rightarrow 0$ in \eqref{proof.lem.biot-savart.4}. Then, recalling the definition of $v_{1,1}$ and \eqref{proof.lem.biot-savart.1}, we take the trace $x_2\rightarrow 0$ and obtain  from \eqref{eq.vorticity.boundary} that
\begin{align}
v_1 (t) & = v_1 (s)   + \int_s^t  \bigg ( \partial_1 p_{F(\tau)} -  \int_{0}^\infty \partial_1  e^{-y_2(-\partial_1^2)^\frac12} (u_1 \omega ) (\tau,\cdot,y_2 ) \dd y_2  \nonumber \\
& ~~~~~~~~~ - \int_{0}^\infty (-\partial_1^2)^\frac12  e^{-y_2(-\partial_1^2)^\frac12} (u_2 \omega ) (\tau,\cdot,y_2 ) \dd y_2 \bigg ) \dd \tau \label{proof.lem.biot-savart.5}
\end{align}
on $\partial\R^2_+$. Since $p_{F(\tau)}$ is the solution given by Proposition \ref{prop.poisson} with $F(\tau) = - u(\tau)\otimes u (\tau)$, by using ${\rm div}\, {\rm div}\, F = - {\rm div}\, (u^\bot \omega ) - \Delta |u|^2/2$ we have the representation
\begin{align}
\partial_1 p_{F(\tau)} = -\partial_1 \frac{|u(\tau )|^2}{2} + \partial_1 \pi_{F(\tau)},\label{proof.lem.biot-savart.6}
\end{align}
where
\begin{align}
\partial_1 \pi_{F(\tau)} & = e^{-x_2 (-\partial_1^2)^\frac12} \int_0^\infty \big ( \partial_1 e^{-y_2 (-\partial_1^2)^\frac12} (u_1 \omega ) (\tau,\cdot,y_2 ) + (-\partial_1^2)^\frac12 e^{-y_2 (-\partial_1^2)^\frac12} (u_2\omega ) (\tau,\cdot,y_2 ) \big ) \dd y_2 \nonumber \\
& ~~+ \int_0^{x_2} \int_{y_2}^\infty  \partial_1 (-\partial_1^2 )^\frac12 e^{-(x_2-2 y_2 + z_2 )  (-\partial_1^2)^\frac12} (u_1 \omega ) (\tau,\cdot,z_2 ) \dd z_2 \dd y_2 \nonumber \\
&  ~~~~- \int_0^{x_2} \int_{y_2}^\infty  \partial_1^2 e^{-(x_2 - 2 y_2 + z_2 ) (-\partial_1^2)^\frac12} (u_2\omega ) (\tau,\cdot,z_2 ) \dd z_2 \dd y_2 \nonumber \\
&  ~~~~~~ -\int_0^{x_2} e^{-(x_2-y_2)(-\partial_1^2)^\frac12} (u_1 \omega ) (\tau,\cdot,y_2 ) \dd y_2.\label{proof.lem.biot-savart.7}
\end{align}
Thus \eqref{proof.lem.biot-savart.5}-\eqref{proof.lem.biot-savart.7} leads to $v_1(t) = v_1 (s)$ on $\partial\R^2_+$ for all $-\infty<s<t<0$. Then \eqref{est.lem.spatial.decay.vorticity} and \eqref{est.lem.temporal.decay.vorticity.1} imply that
\[
v_1(t)=\displaystyle \lim_{s\rightarrow -\infty} v_1 (s)=\lim_{s\rightarrow -\infty} \int_0^\infty e^{-y_2 (-\partial_1^2)^{1/2}} \omega (s,\cdot,y_2) \dd y_2 = 0
\]
on $\partial\R^2_+$ by the Lebesgue convergence theorem. The proof is now complete.

\section{Application to geometric regularity criterion}\label{sec.criterion}

We shall extend a geometric regularity criterion \cite {GiMi} of solutions to the Navier-Stokes equations in $\R^3$ to the case when the domain is the half space $\R^3_+$ with the Dirichlet condition as an application of the Liouville type result (Theorem \ref{thm.intro.1}). As already discussed in \cite {GiMi} when one imposes the Neumann boundary problem (or the slip boundary condition), the extension is rather straightforward. This is because the rescaled two-dimensional vorticity equations still enjoy the maximum principle since there is no vorticity production from the boundary. We shall state our geometric regularity criterion for the Dirichlet problem in a rigorous way.

We consider the Navier-Stokes equations in the half space $\R^3_+ =\{(x_1,x_2,x_3)\in \R^3 \mid x_3 >0\}$
\begin{align}
\partial_t u - \Delta u +\nabla \cdot ( u \otimes u )  + \nabla p =0,~~~~~ {\rm div} ~ u =0~~~~~~~~~{\rm in}~~ (0, T) \times \R^3_+\label{eq.NS.Sec.6}
\end{align}
with the Dirichlet boundary condition:
\begin{align}
u=0~~~~~~~~~~{\rm on} ~~ (0, T) \times \partial \R^3_+.\label{eq.Dirichlet.Sec.6}
\end{align}

As mentioned in the introduction, we need to consider a spatially non-decaying solution to carry out what is called a blow-up argument. However, if one allows non-decaying solutions, the uniqueness of the initial-boundary value problem for \eqref {eq.NS.Sec.6}-\eqref {eq.Dirichlet.Sec.6} fails. Indeed,  the Poiseuille type flow of the form

\begin{align}
u=(u_1 (t,x_3),0,0),~~~p(t,x_1)=-x_{1}f(t),\label{eq.Poiseuille type flow.Sec.6}
\end{align}
solves \eqref {eq.NS.Sec.6}-\eqref {eq.Dirichlet.Sec.6} provided that $u_1$ solves the heat equation
\begin{align*}
\partial_{t}u_1 - \partial_3^2 u_1=f(t) ~~~{\rm in}~~ (0, T) \times \{x_3>0\},\\
u_1=0~~~~{\rm on}~~ (0, T) \times \{x_3=0\}.
\end{align*}
with some $f$ depending only on time. Since $f$ can be chosen arbitrary, one is able to construct various solutions $(u,p)$ to \eqref {eq.NS.Sec.6}-\eqref {eq.Dirichlet.Sec.6} of the form \eqref {eq.Poiseuille type flow.Sec.6} with the same initial data. If one assumes that $f$ is bounded and smooth, all such $(u,p)$ is smooth and bounded. Hence this yields the non-uniqueness of the initial-boundary value problem for \eqref {eq.NS.Sec.6}-\eqref {eq.Dirichlet.Sec.6} when one allows non-decaying solutions.

A simple way to avoid non-uniqueness is to improve a relation between the pressure and the velocity. Taking the divergence of \eqref {eq.NS.Sec.6}, we see 
\begin{align}
-\Delta p=\sum^{3}_{i,j=1}{\partial_i}{\partial_j}(u_i u_j)~~~{\rm in}~~   \R^3_+,  \label{eq.relation between u and p.Sec.6}
\end{align}
since ${\rm div} \, u =0$. Next, taking the inner product of \eqref {eq.NS.Sec.6} with normal $n=(0,0,-1)$, we have
\begin{align}
\frac{\partial p}{\partial n}=-\Delta u\cdot n ~~~{\rm on}~~ \partial  \R^3_+.\label{eq. NS times n.Sec.6}
\end{align}
It is convenient to decompose $p$ into the sum $p_H + p_F$ as we did in earlier sections. Namely, for the harmonic pressure term $p_H$ we require
\begin{align}
-\Delta p_H & =0~~~{\rm in}~~   \R^3_+, \label{eq.harmonic pressure.Sec.6}\\
\frac{\partial p_H}{\partial n} & =-\Delta u\cdot n ~~~{\rm on}~~ \partial  \R^3_+. \label{eq. harmonic pressure boundary.Sec.6}
\end{align}
and the pressure $p_F$ coming from transport term we require
\begin{align}
-\Delta p_F & =\sum^{3}_{i,j=1}{\partial_i}{\partial_j}(u_i u_j),~~F=(u_i u_j)~~~{\rm in}~~   \R^3_+, \label{eq.pressure transport term.Sec.6}\\
\frac{\partial p_F}{\partial n} & =0 ~~~{\rm on}~~ \partial  \R^3_+. \label{eq. pressure transport term.boundary.Sec.6}
\end{align}
Evidently, \eqref {eq.harmonic pressure.Sec.6}-\eqref {eq. pressure transport term.boundary.Sec.6} implies \eqref {eq.relation between u and p.Sec.6}-\eqref {eq. NS times n.Sec.6} for $p=p_H + p_F$. Note that the $\Delta u\cdot n=-{\rm div}_{\partial  \R^3_+}(\omega \times n)$ as noted in \cite {AbGi}. If one imposes  smoothness and boundedness for $u$ up to second derivatives, one can get the uniqueness of $\nabla p$ (determined from $u$) provided that  $p$ is restricted to avoid the linear growth at spatial infinity; see Proposition \ref {prop.poisson} and Proposition \ref {prop.laplace} and also \cite {AbGi}. The unique solution is formally written by using the Helmholtz projection $\mathbb P$ to the solenoidal space:
\begin{align}
\nabla p=(I-\mathbb P)(\Delta u-\nabla \cdot ( u \otimes u ))\label{eq. Helmholtz projection.Sec.6}
\end{align}
and the solution having this form is called a mild solution. It is not difficult to prove the uniqueness of the mild solution; see \cite {GiIMa} for the whole space and \cite {BaJi} for the half space.

There is a large literature giving a growth condition for pressure so that the solution is a mild solution which is unique. Such type of result goes back to \cite{GaMa} and has been developed in the case of the whole space \cite {GIKM} and the half space \cite {Mar08}. A typical criterion for the whole space case is $p\in L^1 ((0,T);BMO(\R^3))$ \cite {Ka}. There are references on this issue \cite{Mar08}, \cite {Mar11} for further relaxation of growth assumptions for the pressure.

In this section we consider the mild solution. We know there is a unique local-in-time mild solution for the initial-boundary value problem for \eqref {eq.NS.Sec.6}-\eqref {eq.Dirichlet.Sec.6} for any bounded continuous initial velocity $u_0$ i.e., $u_0 \in BC(\overline{\R^3_+})$ which is solenoidal in the sense that ${\rm div} \, u_0 =0$ in $\R^3_+$ and $u_0 \cdot n=0$ on $\partial \R^3_+$ \cite {So03}; see also \cite {BaJi}.

We are now in position to prove Theorem \ref{thm.intro.2}, which is a natural extension of the geometric regularity criterion of \cite{GiMi}. 
We shall prove this result by a blow-up argument. The basic strategy is the same as in \cite {GiMi}. However, to assert uniqueness of the limit we invoke our Liouville type result (Theorem \ref {thm.intro.1}).
Of course, in some steps it is more involved because of the presence of the boundary.

\noindent {\it Proof of Theorem \ref{thm.intro.2}.} Step 1 (Construction of blow-up sequence). Assume that $u$ blows up at $t=T$. Then there exists a sequence $\{(t_k,x_k)\}^{\infty}_{k=1}\subset [0,T)\times \R^3_+$ with $t_{k+1}>t_k$ such that\\
(i) $|u(t,x)|\leq M_k$ for $t \leq t_k$,~~$x \in \R^3_+$\\
(ii) $M_k=\|u (t_k) \|_{\infty}\rightarrow \infty$,~~$t_k \uparrow T$ as $k\rightarrow \infty$\\
(iii) $|u(t_k,x_k)|\geq M_k /2$\\
We rescale $u$, $\omega$ with respect to $(t_k,x_k)$ i.e. 
\begin{align*}
u_k (t,x)={\lambda}_k u(t_k +{{\lambda}_k}^2 t, x_k + {\lambda}_k x)~~~~~~~~~~~~~~~~~~~~~~~~~~~~~~\\
{\omega}_k (t,x)={{\lambda}_k}^2 \omega (t_k +{{\lambda}_k}^2 t, x_k + {\lambda}_k x),~~T-t_k >{{\lambda}_k}^2 t>-t_k 
\end{align*}
with ${\lambda}_k=1/M_k$. Since \eqref {eq.NS.Sec.6}-\eqref {eq.Dirichlet.Sec.6} is scaling invariant under the above rescaling, we see that $u_k$ is a mild solution of \eqref {eq.NS.Sec.6}-\eqref {eq.Dirichlet.Sec.6} in $(-t_k {M_k}^2,0]\times \R^3_{+,-c_k}$ with $c_k=x_{k,3} M_k$, where $x_k=(x_{k,1},x_{k,2},x_{k,3})$ and $\R^3_{+,-c}=\{x=(x_1,x_2,x_3)\in \R^3 \mid x_3 >-c\}$.

\noindent{Step 2 (Compactness).} By assumption (i) we have $|u_k|\leq 1$ in $(-t_k {M_k}^2,0]\times \R^3_{+,-c_k}$. Since $u_k$ is a mild solution, we know that $\nabla u_k$ is also bounded in $(-t_k {M_k}^2 +1,0]\times \R^3_{+,-c_k}$ by a result of \cite {BaJi}. Thus $(u_k, {\omega}_k)\rightharpoondown (\bar{u},\bar{\omega})$ as $k\rightarrow \infty$ *-weakly in $L^{\infty}$ with some $(\bar{u},\bar{\omega})$ such that $|\bar{u}|\leq 1$ and $|\bar{\omega}|\leq c$ in $(-\infty ,0]\times \R^3_{+,-c}$ $(c=\lim_{k\rightarrow \infty}c_k)$ by taking a subsequence. Moreover, $\bar u$ is a bounded global mild solution in $(-\infty ,0]\times \R^3_{+,-c}$.
Note that there are two cases depending upon whether $\lim c_k=\infty$ or $\lim c_k <\infty$. In the first case $\R^3_{+,-c}=\R^3$ and the limit $\bar u$ solves the Navier-Stokes equations in the whole space. In the second case $\bar u$ solves the Navier-Stokes equations in the half space $\R^3_{+,-c}$ with the Dirichlet condition (cf. \cite{Gi1}).

We need some compactness to guarantee that $u_k$ converges to $u$ at least locally uniformly in $(-\infty ,0]\times \R^3_{+,-c}$ to guarantee that $u_k (0,0)\rightarrow {\bar u}(0,0)$.

In the whole space this can be guaranteed by the estimates of higher-order derivatives so that all space-time derivatives of  $u_k$ are bounded in $(-t_k {M_k}^2 +1,0]$ uniformly in $k$ (e.g. \cite {GiSa}). In the case of the Dirichlet problem it seems  to be unknown since it is nontrivial to handle normal derivatives. However, what we need here are local estimates, rather than global estimates.

We first note that the pressure defined by \eqref {eq. Helmholtz projection.Sec.6} is estimated as
\begin{align}
\|p\|_{L^r (B_R (x_0)\cap \R^3_+)}\leq C(\|\omega \|_{L^{\infty}(\partial \R^3_+)}+\|u\|^2_{L^\infty (\R^3_+)}) \label{eq. pressure estimate.Sec6}
\end{align}
with $C$ depending on $R$ and $r\in (1,\infty)$ and independent of $u$ and $\omega$, where $B_R (x_0)$ is a closed ball of radius $R$ centered at $x_0\in \R^3_+$. Here we normalize $p$ such that $p(x_0)=0$. Decompose $p$ into $p_H+p_F$. For $p_F$ we have a BMO estimate $\|p\|_{BMO} \leq C\|u\|^2_{\infty}$. For the harmonic pressure term, as observed in \cite {AbGi}, we have 
\begin{align*}
\|x_3 \nabla p_H\|_{L^{\infty}(\R^3_+)}\leq C\|\omega\|_{L^{\infty}(\partial \R^3_+)}.
\end{align*}
From these two estimates \eqref {eq. pressure estimate.Sec6} easily follows. The estimate \eqref {eq. pressure estimate.Sec6} enables us to localize the problem. We cut off $u$ in $B_R (x_0)\cap \R^3_+$ with Bogovski type adjustment to apply the $L^r$ maximal regularity of the Stokes equation problem in a smoothly bounded domain with the zero boundary condition, e.g. \cite {GiSo}. By \eqref {eq. pressure estimate.Sec6} we observe that the external force has a local space-time $L^r$ bound depending on $u$ only through the space time sup norm of $u$ and $\nabla u$. Thus we are able to control all $W^{2,1}_r (I \times (B_R (x_0)\cap \R^3_+))$ norm of $u$, where $I$ is a bounded time interval $\subset (-\infty,0]$. By the Sobolev embedding theorem we have a H\"{o}lder bound on $\nabla u$ in $Q=I \times (B_R (x_0)\cap \R^3_+)$. This is of course enough to ensure that $u_k$ converges to $\bar u$ locally uniformly in $(-\infty ,0]\times {\bar{\R}}^3_{+,-c}$. By a bootstrap argument we improve the regularity of the pressure and observe that $u_k\rightarrow u$ locally uniformly for its all derivatives. Note that without a bound for the pressure one cannot localize the problem.
Since $(t_k,x_k)$ is taken so that $|u_k (0,0)|\geq 1/2$ by Step 1 (iii), we conclude that $|{\bar u}(0,0)|\geq 1/2$.

\noindent{Step 3 (Characterization of the limit).} We now apply the continuous alignment condition {\rm (CA)} and our Liouville type result (Theorem \ref {thm.intro.1}) to conclude that $\bar u$ must be zero, which contradict with $|{\bar u}(0,0)|\geq 1/2$.
Here is a sketch of the proof. We set the vorticity direction $\xi_k=\omega_k /|\omega_k|$. Then {\rm (CA)} implies
\begin{align*}
|\xi_k (t,x)-\xi_k (t,y)|\leq \eta (\frac{|x-y|}{M_k})\rightarrow 0,
\end{align*}
so that $\bar{\xi}=\bar{\omega} / |\bar{\omega}|$ is independent of $x$. By the unique existence theory \cite {BaJi} of the mild solution $\bar{\xi}$ must be also constant in time. Thus $(\bar u, \bar{\omega})$ is a two-dimensional flow in $(-\infty ,0)\times \R^3_{+,-c}$. When $c=\infty$ the problem is reduced to the whole space case, and it is already proved  in \cite{GiMi} that $\bar{u}=0$, which leads to a desired contradiction. Hence it suffices to consider the case $c<\infty$. By a suitable change of coordinates we may assume that $\bar u = (\bar{u}_1 (x_1,x_2),\bar{u}_2 (x_1,x_2),0)$ with $\bar{\omega} = (0,0,\bar{\omega}_3)$, $\bar{\omega}_3 \geq 0$ and $\bar{u}_1=\bar{u}_2=0$ on $(-\infty ,0)\times \partial \R^2_{+}$ where $\R^2_+ =\{(x_1,x_2)\in \R^2 \mid x_2 >0\}$. 

Now we shall apply Theorem \ref{thm.intro.1} for $(\bar u,\bar{\omega})$. The condition {\rm (C4)} is trivially fulfilled because $\bar{\omega_3}\geq 0$. The condition {\rm (C3)} is inherited from the type I assumption. It remains to prove {\rm (C1)} and {\rm (C2)} for our mild solution $\bar u$, but thanks to Proposition \ref{prop.poisson} and \ref{prop.laplace}, it is enough to prove {\rm (C1)}. 
By the construction of the blow-up sequence we know 
\begin{align*}
\sup_{-\infty <t<0}\|\bar u (t) \|_{\infty} <\infty.
\end{align*}
Applying a result of \cite{BaJi}, we also know $\|\nabla \bar u (t) \|_{\infty}$ is bounded for all $t<0$. Then we have to estimate the higher-order derivatives to prove {\rm (C1)}, which will be established in  Lemma \ref{lem.hsu} below. This is  sufficient to derive {\rm (C1)} so we apply Theorem \ref{thm.intro.1} to conclude $\bar u \equiv 0$, and reach a contradiction.

\begin{lem}\label{lem.hsu} Let $u$ be a mild solution of \eqref {eq.NS.Sec.6}-\eqref {eq.Dirichlet.Sec.6} in $(0,T)\times \R^2_+$ with initial data $u_0 = (u_{0,1}, u_{0,2})$ in $BC(\overline{\R^2_+})$. Assume that $u_0$ is solenoidal, i.e. ${\rm div}\, u_0 =0$ in $\R^2_+$ and $u_{0,2} =0$ on $\partial \R^2_+$. Assume that there exists $T>0$ such that
\begin{align*}
K=\sup_{0<t<T}\|u (t) \|_{\infty} <\infty.
\end{align*}
Then there exists a constant $C$ depending only on $K$ and $T$ such that 
\begin{align*}
\sup_{0<t<T}(t^{\frac{m}{2}}\|{\nabla}^m u  (t) \|_{\infty} + t^{1+\frac{l}{2}}\|{\nabla}^l \partial_t u (t) \|_{\infty})\leq C
\end{align*}
with $m=1,2,3$ and $l=0,1$, where ${\nabla}^0$ is interpreters as an identity operator.
\end{lem}

\begin{rem}{\rm We need the H\"{o}lder norm estimates in (C1), but these are obtained from a simple interpolation of $L^{\infty}$ bounds in Lemma \ref{lem.hsu}; see e.g. \cite[Theorem 3.2.1]{Kry}, \cite[Section 3.2]{Ta}.
}
\end{rem}

The idea of the proof of Lemma \ref{lem.hsu} is to estimate the tangential derivatives with up to one normal derivative as in \cite {GiGiSa} or \cite {GiSa}. We also need to estimate the time derivative. In the meanwhile we estimate $p_H$ and $p_F$, which enable us to estimate the normal derivatives. Except for the estimates of the pressure term the argument is rather conventional, so we give a sketch of the proof instead of giving a full detail. In the argument below  we just use $L^{\infty}$ norm so we simply write $\|f\|$ instead of $\|f\|_{\infty}$.

\noindent {\it Sketch of the proof of Lemma \ref{lem.hsu}.} Step 1 (Tangential derivatives and time derivatives). We first note that the mild solution solves the integral equation
\begin{align}\label{stokes semigroup}
u(t)=S(t)u_{0}+w,~~~~ w=-\int_{s}^{t} S(t-s)\mathbb{P}\nabla \cdot ( u \otimes u ) (s) \dd s,
\end{align}
where $S(t)$ is the Stokes semigroup. According to \cite {BaJi}, we know
\begin{align}\label{BaeJin est}
\|\nabla S(t)\mathbb{P}f\|\leq C_1 t^{-\frac12}\|f\|,~~~~ ~\|S(t)\mathbb{P}\nabla \cdot f\|\leq C_2 t^{-\frac12}\|f\|.
\end{align}
Taking the derivatives $\nabla$ in \eqref{stokes semigroup}, we obtain, for $0<\epsilon<1$, 
\begin{align*}
\|\nabla u\| (t)   & \leq  C_{0}t^{-\frac12}\|u_0\| +\int_{t(1-\epsilon)}^{t}\| \nabla S(t-s)\mathbb{P} \nabla \cdot ( u \otimes u ) (s) \| \dd s \\
& ~~+ \int_{0}^{t(1-\epsilon)}\| \nabla S(\frac{t-s}{2})\cdot S(\frac{t-s}{2}) \mathbb{P}\nabla \cdot ( u \otimes u ) (s) \| \dd s\\
& \leq  C_{0}t^{-\frac12}\|u_0\|+2KC_{1}\int_{t(1-\epsilon)}^{t}(t-s)^{-\frac12}\| \nabla u (s) \|  \dd s +K^{2} C_1 C_2 \int_{0}^{t(1-\epsilon)}(t-s)^{-1} \dd s.
\end{align*}
This yields the estimate
$\|\nabla u (t) \| \leq  C_{K}t^{-1/2}$, $t\in (0,T)$
by using \cite[Lemma 2.4]{GiGiSa}.
Since the tangential derivative $\partial_1$ commutes with the Stokes semigroup, a similar argument yields $
\|\partial_{1}^{m-1} \nabla u (t) \| \leq  C_{Km}t^{-m/2}$, $t\in (0,T)
$ for all $m=1,2,\cdots$. Note that the proof makes the sense if we know in advance that
$\| \partial_1^{m-1} \nabla u (t) \|$ is finite and locally bounded in time in $(0,T)$; 
however we are able to justify this process by approximating $u_0$ by $L^{p}_{\sigma}$ vector field, and the details are omitted here.

As for the time derivative, we differentiate \eqref{stokes semigroup} in $t$, which gives
\begin{align*}
\partial_{t}u (t) & =\frac{\dd S(t)}{\dd t}u_{0}-\int_{t(1-\epsilon)}^{t}  S(t-s)\mathbb{P}\partial_{s}\nabla \cdot ( u \otimes u ) (s) \dd s -S(\epsilon t)\mathbb{P}\nabla \cdot ( u \otimes u )(t-\epsilon t)\\
& ~~~ - \int_{0}^{t(1-\epsilon)}\frac{\dd}{\dd t} S(t-s) \mathbb{P}\nabla \cdot ( u \otimes u ) (s) \dd s.
\end{align*}
By using the estimate \eqref{BaeJin est} and $\| \frac{\dd S(t)}{\dd t}f\|\leq C t^{-1} \|f \|$ that is obtained from the explicit formula of the Stokes semigroup in \cite{So03,U},  we have
\begin{align*}
\|\partial_{t}u (t) \| & \leq C_{K}t^{-1}+\int_{t(1-\epsilon)}^{t}\frac{C_K}{(t-s)^{\frac12}}\|\partial_{s}u (s) \|\dd s 
+C_{K}(\epsilon t)^{-\frac12}+\int_{0}^{t(1-\epsilon)}\frac{C_K}{(t-s)^{\frac32}}\dd s\\
& \leq C_{K,T,\epsilon}t^{-1}+ \int_{t(1-\epsilon)}^{t}\frac{C_{K}}{(t-s)^{\frac12}} \|\partial_{s}u (s) \|\dd s. 
\end{align*}
Then by using \cite[Lemma 2.4]{GiGiSa} we have the estimate $\|\partial_{t}u (t) \|\leq C_{K,T} t^{-1}$. Similarly, we can also obtain the following estimate.
\begin{align*}
\|\nabla \partial_{t}u (t) \| & \leq \|\nabla \frac{\dd S(t)}{\dd t}u_0\|+\int_{t(1-\epsilon )}^{t}\|\nabla S(t-s)\mathbb{P}\partial_{s}\nabla \cdot ( u \otimes u ) (s) \|\dd s \\
& ~~~ +\|\nabla S(\epsilon t)\mathbb{P}\nabla \cdot ( u \otimes u )(t-\epsilon t)\|+\int_{0}^{t(1-\epsilon)}\|\nabla \frac{\dd}{\dd t} S(t-s) \mathbb{P}\nabla \cdot ( u \otimes u )(s) \|\dd s.
\end{align*}
By using \eqref{BaeJin est} again, we have 
\begin{align*}\|\nabla \partial_{t}u (t) \|\leq C_{K,T,\epsilon}t^{-\frac32}+\int_{t(1-\epsilon)}^{t}\frac{C_{K }}{(t-s)^{\frac12}} \|\nabla \partial_{s}u (s) \|\dd s.
\end{align*}
Hence, by \cite[Lemma 2.4]{GiGiSa} we arrive at  $\|\nabla \partial_{t}u (t) \|\leq C_{K,T} t^{-3/2}$.

\noindent Step 2 (Pressure estimates). In order to estimate the normal derivatives of the solution, we first consider $p_F$ and $p_H$. Recall that $p_F$ is expressed as $p_F=\int_{\R^2}\nabla \nabla \cdot E(x-y) ( \tilde u \otimes \tilde u )\dd y$; see the proof of Proposition \ref{prop.poisson}. Hence $\partial_{1}p_F$ can be decomposed into $\int_{R^2} \nabla^3 E(x-y)( \tilde u \otimes \tilde u ) (1-\chi_R) \dd y +\int_{\R^2} \nabla E(x-y) \partial_{1}\nabla \cdot \big ( \tilde u \otimes \tilde u   \chi_R \big ) \dd y$, where $E(x)$ is the newton potential and $\chi_R=\chi _R(x-y)$ is a smooth cut-off such that $\chi_R (x-y)=1$ for $|x-y|\leq R$ and $\chi_R (x-y) =0$ for $|x-y|\geq 2 R$. Then we have 
\begin{align*}
\|\partial_1 p_F (t)  \|\leq R^{-1} \|u (t) \| ^2 +R\big  (\|\partial_{1}\nabla u (t) \| \|u (t) \| + \|\nabla u (t) \|^2 \big ),
\end{align*}
which yields $\|\partial_{1}p_F \|\leq Ct^{-1/2}$ by taking $R=t^{1/2}$, where $C$ depends only on $K$ and  $T$. As for $p_H$, we see
\begin{align*}
\|\partial_{1}p_H  \|_{L^\infty_{x_1}} (t, x_2) \leq \int_{x_2}^{L}\|\partial_{1}\partial_{2} p_H \|_{L^\infty_{x_1}}  (t, y_2)  \dd y_2 + \|\partial_{1} p_H \|_{L^\infty_{x_1}} (t, L) .
\end{align*}
Since $\partial_{1}\partial_{2} p_H (w) =\partial_{1}^{2} e^{- x_2 (-\partial_1^2)^{1/2}} (\omega |_{x_2=0} ) = e^{- x_2 (-\partial_1^2)^{1/2}} (\partial_{1}^{2} \omega |_{x_2=0})$, we have
\begin{align}\label{harmonic pressure estimate}
\|\partial_{1}\partial_{2} p_H\|_{L^\infty_{x_1}} (t, y_2)\leq C \| \partial_{1}^{2}\nabla u (t) \|.
\end{align}
Furthermore,  it follows from \eqref{est.prop.laplace.1} that $\|\partial_{1} p_H\|_{L^\infty_{x_1}} (t, L)\leq  C L^{-1} \|\omega (t) \|$. Hence 
\begin{align*}
\|\partial_{1} p_H\|_{L^\infty_{x_1}} (t, x_2)\leq C \big (L\| \partial_{1}^{2}\nabla u (t) \|+\frac{1}{L}\|\omega (t) \| \big ), ~~~~ 0\leq x_2 \leq L.
\end{align*}
By taking $L=t^{1/2}$, we have $\sup_{0\leq x_2 \leq L}\|\partial_{1} p_H \|_{L^\infty_{x_1}} (t, x_2) \leq Ct^{-1}$. Then by maximum principle we obtain $\|\partial_{1} p_H (t) \|\leq Ct^{-1}$. The constant $C$ depends only on $K$ and $T$. From the similar argument we can extend the estimates to the higher-order tangential derivatives.

\noindent Step 3 (Normal derivatives). By combining the above estimates with the following equation
\begin{align}\label{first eq}
\partial_{t}u_1 -\Delta u_1 + (u,\nabla)u_1 +\partial_{1}p=0,
\end{align}
it is easy to check that $\|\partial_{2}^{2}u (t) \|\leq Ct^{-1}$. Finally, by differentiating \eqref{first eq} in the normal direction and by using \eqref{harmonic pressure estimate}, the estimate $\|\partial_{2}^{3}u_1 (t) \|\leq Ct^{-3/2}$ follows. 
With the aid of the estimate $\|\partial_{1}^{2}p (t) \|\leq Ct^{-3/2}$ and the divergence free property of the solutions, we finally obtain $\|\partial_{2}^{3}u_2 (t) \|\leq Ct^{-3/2}$ by differentiating \eqref{first eq} in the tangential direction. The proof is now complete.

\end{document}